\documentclass[10pt,twoside, a4paper, english, reqno]{amsart}
\usepackage{graphicx, amsmath, varioref, amscd, amssymb,color, bm, stmaryrd, amsthm}
\usepackage{fancybox}
\usepackage[pdftex, colorlinks=false, backref]{hyperref}

\numberwithin{equation}{section}
\allowdisplaybreaks

\newtheorem{theorem}{Theorem}[section]
\newtheorem{lemma}[theorem]{Lemma}
\newtheorem{corollary}[theorem]{Corollary}

\newtheorem{proposition}[theorem]{Proposition}
\theoremstyle{definition}
\newtheorem{definition}[theorem]{Definition}

\theoremstyle{remark}

\newcommand{\Curl}{\operatorname{curl}}
\newcommand{\Grad}{\nabla}
\newcommand{\vr}{\varrho}
\newcommand{\vu}{\vc{u}}

\newcommand{\pvint}{{\rm P.V.}\!\int}

\newcommand{\vc}[1]{{\bm{#1}}}

\newcommand{\dott}{\,\cdot\,}

\newcommand{\norm}[1]{\left\Vert#1\right\Vert}
\newcommand{\abs}[1]{\left|#1\right|}

\newcommand{\Dt}{\Delta t}
\newcommand{\R}{\mathbb{R}}
\newcommand{\N}{\mathbb{N}}

\newcommand{\Om}{\ensuremath{\Omega}}

\newcommand{\vt}{\theta}
\newcommand{\vvt}{\vartheta}

\begin{document}
\title[Operator splitting for 2D incompressible fluid equations]
{{Operator splitting for two-dimensional incompressible fluid equations}}

\author[Holden]{Helge Holden}
\address[Holden]{\newline
    Department of Mathematical Sciences,
    Norwegian University of Science and Technology,
    NO--7491 Trondheim, Norway,\newline
{\rm and} \newline
  Centre of Mathematics for Applications, 
University of Oslo,
  P.O.\ Box 1053, Blindern,
  NO--0316 Oslo, Norway }
\email[]{\href{holden@math.ntnu.no}{holden@math.ntnu.no}}
\urladdr{\href{http://www.math.ntnu.no/~holden}{www.math.ntnu.no/\~{}holden}}

\author[Karlsen]{Kenneth H.~Karlsen}
\address[Karlsen]{\newline
   Centre of Mathematics for Applications, 
 University of Oslo,
  P.O.\ Box 1053, Blindern,
  NO--0316 Oslo, Norway}
\email[]{kennethk@math.uio.no}
\urladdr{www.math.uio.no/\~{}kennethk}

\author[Karper]{Trygve Karper}
\address[Karper]{\newline
    Department of Mathematical Sciences,
    Norwegian University of Science and Technology,
    NO--7491 Trondheim, Norway}
\email[]{\href{karper@math.ntnu.no}{karper@math.ntnu.no}}
\urladdr{\href{http://www.math.ntnu.no/~karper}{www.math.ntnu.no/\~{}karper}}

\date{\today}




\thanks{Supported in part by the Research Council of Norway.}

\maketitle

\begin{abstract}
We analyze splitting algorithms for a class of 
two-dimensional fluid equations, which includes the incompressible 
Navier--Stokes equations and the surface quasi-geostrophic equation. 
Our main result is that the Godunov and Strang 
splitting methods converge with the expected rates
 provided the initial data are sufficiently regular. 
\end{abstract}

\section{Introduction}
Let $T > 0$ be a finite final time. We are interested in solutions 
$\vt\colon [0,T]\times \R^2\to \R$ to the generalized active scalar equation
\begin{equation}\label{eq:qg}
	\vt_t + \vu \cdot \Grad \vt 
	+ \Lambda^\alpha \vt = 0, 
	\quad \text{in }(0,T)\times \R^2,\quad \alpha \in (0,2],
\end{equation}
where  $\Grad$ is the  gradient operator and $\Lambda=(-\Delta)^{1/2}$ is the fractional Laplacian defined through 
Riesz operators (see Section~\ref{sec:2}).
The divergence free  velocity $\vu$ is determined directly from $\vt$ through the nonlocal relation
\begin{equation*}
	\vu = \Curl \Lambda^{-\beta}\vt, \qquad 
	\beta \in [1,2],
\end{equation*}
where $\Curl=\nabla_x^\bot$ denotes the spatial curl operator defined
by $\Curl(\phi)=(-\phi_y, \phi_x)$.


The general active scalar equation \eqref{eq:qg} seems  
to have appeared first in the mathematical literature in \cite{Constantin}. The 
general formulation encompasses a whole class of two-dimensional 
fluid equations, interpolating between the Euler/Navier--Stokes 
equations and the surface quasi-geostrophic equation. 
Different choices of $\alpha$ and 
$\beta$ lead to different fluid equations. 
The most interesting (and studied) examples are the \emph{Navier--Stokes}
equations ($\alpha= \beta = 2$) and the \emph{surface quasi-geostrophic}
equation ($\beta = 1$,  $0< \alpha < 2$) \cite{Constantin:1994uq};
in the latter case, $\vu=(-\mathcal{R}_2\theta,\mathcal{R}_1\theta)$ with 
$\mathcal{R}_1,\mathcal{R}_2$ denoting the usual Riesz transforms in $\R^2$. 

The quasi-geostrophic equation has recently 
received considerable attention from the mathematical community. 
In particular, since the two-dimensional Navier--Stokes 
equations admit smooth solutions, it has been a question 
if the quasi-geostrophic equation exhibits similar behavior. 
In this respect, it is common to distinguish between 
three cases of $\alpha$ for the geostrophic equation. 
When $\alpha \in (1,2)$, the dissipation term $\Lambda^\alpha$ provides enough 
regularization to guarantee the existence of smooth solutions \cite{Constantin2}. 
When $\alpha \in (0,1)$, the global properties 
of solutions are still open. The remaining case ($\alpha = 1$) is 
known as the critical case, and the global behavior of solutions was 
settled only recently (cf.~\cite{Caffarelli1,Kiselev:2007fk,Dong} and the references therein). 

In terms of the physical applicability of \eqref{eq:qg}, 
it seems like the most relevant models are the critical geostrophic equation and 
of course the Navier--Stokes equations.  In particular, the critical geostrophic 
equation has been proposed as a simplified model
for strongly rotating atmospheric flow. We refer 
the reader to \cite{Pedlosky:1987fk} for more on the 
physical aspects of the model.

We now turn to the main topic of the present paper, namely operator 
splitting algorithms for computing approximate solutions to \eqref{eq:qg}. 
Generally speaking, the label ``operator splitting" alludes 
to the well-known idea of constructing numerical methods for an intricate partial 
differential equation by reducing the original equation to a succession of
simpler equations, each of which can be handled by some efficient 
and tailor-made numerical method.  The operator splitting 
approach has been comprehensively described in 
a large  number of articles and books. 
We do not survey the literature here, referring the 
reader to the bibliography in \cite{Holden:book}.

Regarding the Navier--Stokes equations, which is a special case of \eqref{eq:qg}, 
operator splitting (viscous splitting) has been analyzed 
and applied in a great number of works, see, e.g., the 
book by Majda and Bertozzi \cite{Majda:2002kx}. Indeed, viscous splitting has 
been frequently utilized as a design principle for numerical methods for 
the Navier--Stokes equations, including vortex or particle methods and 
transport-diffusion or characteristic-Galerkin methods.  
Error estimates for Godunov and Strang viscous splitting 
algorithms have been established, e.g., in  \cite{Majda:2002kx},  utilizing 
arguments that are different from ours and which rely on the well-known 
fact that there exist unique, global smooth solutions to the two dimensional Euler and 
Navier--Stokes equations.

In this paper we apply operator splitting to separate the effects 
in \eqref{eq:qg}  of the transport term $\vu = \Curl \Lambda^{-\beta}\vt\cdot \Grad \vt$ and 
the fractional diffusion term $\Lambda^\alpha \vt$. 
This type of splitting is reasonable 
as it allows for specialized hyperbolic methods to be applied 
in the transport step and specialized ``Fourier space" 
methods in the fractional diffusion step. The interested reader 
can consult \cite{Holden:book} for further information on operator splitting. 

Our main contribution is that we contribute rigorous proofs of the
expected convergence rates for operator splitting 
applied to the general active scalar equation \eqref{eq:qg}.
Our approach is inspired by the recent paper \cite{Holden:tao} (see also \cite{lub}) 
on splitting algorithms for the KdV equation, and for that reason our results apply 
under the standing assumption that there exists 
a smooth solution to \eqref{eq:qg}. This assumption is verified for the 
Navier--Stokes equations and the quasi-geostrophic equation 
with $\alpha \geq 1$. It is  also reasonable to expect the existence of unique smooth 
solution in the regime $\alpha \in [1,2]$ and $\beta \in [1,2]$, cf.~\cite{Constantin:2008uq,Kiselev:2010vn,Kiselev:2010ys,Miao:2009kx} for results in that direction.

Let us now discuss our splitting methods in more details. For this purpose, 
we write \eqref{eq:qg} in the form:
\begin{equation*}
	\vt_t = C(\vt), \qquad C(\vt) = A(\vt) + B(\vt), 
\end{equation*}
\begin{equation*}
	B(\vt) = -\Curl \Lambda^{-\beta}\vt \cdot \Grad \vt, \qquad A(\vt) = - \Lambda^{\alpha}\vt.
\end{equation*}
We will need the solution operators $\Phi_A(t,  \vt_0)$ and $\Phi_B(t,  \vt_0)$, defined as 
the solutions to the abstract differential equations:
\begin{align*}
	\partial_t\Phi_A(t, \vt_0) = A(\Phi_A(t, \vt_0)), \qquad \Phi_A(0, \vt_0) = \vt_0,\\
	\partial_t\Phi_B(t, \vt_0) = B(\Phi_B(t, \vt_0)), \qquad \Phi_A(0, \vt_0) = \vt_0.
\end{align*}
For local-in-time existence results for the inviscid equation 
$\vt_t + \vu \cdot \Grad \vt =0$ (i.e., existence of the $\Phi_B$ operator) 
when the initial data belong to Sobolev or 
Triebel-Lizorkin spaces and $\beta=1$, see \cite{Constantin:1994uq,Cordoba,Chae:2003fk}; 
although for $\beta\neq 1$ such results cannot be found in the literature, they can be proved by 
properly adapting the arguments in  \cite{Constantin:1994uq,Cordoba,Chae:2003fk}. 
Regarding the $\Phi_A$ operator, the fractional diffusion equation
$$
v_t + \Lambda^\alpha \vt = 0, \qquad v(0)=v_0,
$$
has a solution given by $v(t)=\mathcal{G}_\alpha(t)\star v_0$, where $\mathcal{G}_\alpha(t,x)$ is 
the fundamental solution in $\R^2$ which can be expressed in terms of the Fourier transform
$\widehat{\mathcal{G}}_\alpha(t,\xi)=e^{-t |\xi|^\alpha}$. If $G_\alpha(x)$ denotes 
the inverse Fourier transform of $e^{- |\xi|^\alpha}$, then 
$$
\mathcal{G}_\alpha(t,x)=t^{-\frac{2}{\alpha}}
G_\alpha\left(xt^{-\frac{1}{\alpha}}\right).
$$
We refer to \cite{Landkof:1972fk,Stein:1970pr} for further information.

The first operator splitting method we will study is  known in the literature as 
\textit{Godunov splitting}. The method is defined as follows: 
Set $\vt^0 = \vt_0$ and sequentially determine approximate solutions 
$\vt^n$, $ n=1, \ldots, M$, satisfying
\begin{equation*}
	\vt^n = \Phi_A(\Delta t,\Phi_B(\Delta t,  \vt^{n-1})) 
	= \Phi_A(\Delta t) \circ\Phi_B(\Delta t)(\vt^{n-1}).
\end{equation*}
Formally, one can show that $\|\vt(n\Delta t) - \vt^n \| = \mathcal{O}(\Delta t)$
in an appropriate spatial norm in the limit $\Delta t\to0$ and $n\Delta t=t$. 
In Section \ref{sec:4}, we will rigorously prove 
this linear convergence rate. Specifically, we show that (for small $\Delta t$)
\begin{equation*}
	\|\vt(n\Delta t) - \vt^n \|_{H^{k-2}} \leq C(T)\Delta t\|\vt_0\|_{H^k}, 
 \end{equation*}	
for all $k \geq 5$.
The second method we will consider is known as \textit{Strang splitting}: 
Let $\vt^0 = \vt_0$ and determine sequentially
\begin{equation*}
	\begin{split}
	\vt^{n+1} &= 
	\Phi_B\Big(\frac{\Delta t}{2}, \Phi_A\big(\Delta t, \Phi_B(\frac{\Delta t}{2}, \vt^{n})\big)\Big) \\
	&= \left[\Phi_B(\frac{\Delta t}{2})\circ \Phi_A(\Delta t)\circ \Phi_B(\frac{\Delta t}{2})\right](\vt^n).
	\end{split}
\end{equation*}
Formally, one can show that the method is second order 
in $\Delta t$. That is, $\|\vt(n\Delta t) - \vt^n \| = \mathcal{O}(\Delta t^2)$
in an appropriate spatial norm in the limit $\Delta t\to0$ 
and $n\Delta t=t$.  In Section \ref{sec:4}, we show that 
\begin{equation*}
	\|\vt(n\Delta t) - \vt^n \|_{H^{k-3\alpha}} 
	\leq C\left(\|\vt_0\|_{H^k}\right) \Delta t^2, 
\end{equation*}
for any sufficiently small time step $\Delta t >0$ and for 
all $k \geq \max\{5, 3\alpha\}$.

To build fully discrete numerical methods for the 
fluid equation \eqref{eq:qg}, we have 
to replace the exact solutions operators $\Phi_A$ and $\Phi_B$ by 
appropriate numerical methods.  However, we will not discuss that here.

The paper is organized as follows: 
Section \ref{sec:2} is of an introductory nature and 
collect some results to be used later on.
The convergence rate result for the Godunov splitting is 
proved in Section \ref{sec:godunov}, while the Strang splitting 
is analyzed in Section \ref{sec:strang}.

\section{Preliminary material} \label{sec:2}

\subsection{Existence and regularity results}
Presently there is no complete existence theory 
for the equation \eqref{eq:qg}. In this paper, we will  
assume the existence of a unique solution 
with the same regularity as the initial data. 
In the literature, one can find results confirming this 
assumption for some specific cases of $\alpha$ and $\beta$.

For $\alpha= \beta = 2$, the equation \eqref{eq:qg} is 
the incompressible Navier--Stokes equations. The following result
is by now classical (cf.~\cite{Lady}). 
\begin{theorem}[Navier--Stokes]\label{lem:}
Let $\alpha = \beta = 2$ and $T>0$.
If $\vt_0 \in H^k$, there exists a unique solution 
$\vt \in C(0,T;H^k)$ of \eqref{eq:qg} with initial data $\vt|_{t=0}=\vt_0$.
\end{theorem}

The next theorem gives 
the well-posedness of the quasi-geostrophic equation 
when $\alpha \in [1,2)$. The sub-critical case ($\alpha > 1$) was  
established in \cite{Constantin2}. Well-posedness in the 
critical case $\alpha=1$ can be found in \cite{Caffarelli1,Kiselev:2007fk,Dong}.

\begin{theorem}[quasi-geostrophic]
Let $\beta = 1$, $\alpha \in [1,2)$, and let $T>0$ be a final time. 
Assume that $\vt_0 \in H^k$ for some $k \geq 1$. Then, 
there exists a unique solution $\vt \in C(0,T;H^k)$ 
of \eqref{eq:qg} with initial data $\vt|_{t=0}=\vt_0$.
\end{theorem}

Since \eqref{eq:qg} appears to be better behaved when $\beta > 1$, it is reasonable
to expect that the previous theorem continue to hold in the entire range $\alpha,\beta \in [1,2]$, 
cf.~\cite{Constantin:2008uq,Kiselev:2010vn,Kiselev:2010ys,Miao:2009kx} 
for some relevant results.

In what follows we shall also need a local-in-time existence results 
for Sobolev regular solutions to the inviscid version of \eqref{eq:qg}. 
However, as mentioned in the introduction, the currently available 
results  \cite{Constantin:1994uq,Cordoba,Chae:2003fk}
apply only to the inviscid quasi-geostrophic equation ($\beta=1$)
$$
\vt_t + \vu \cdot \Grad \vt =0,  \qquad 
\vu=(-\mathcal{R}_2\theta,\mathcal{R}_1\theta),
$$
with initial data $\theta|_{t=0}=\theta_0$ belonging to some Sobolev space $H^k$. 
Throughout this paper we will simply make the standing assumption that the 
inviscid generalized quasi-geostrophic 
equation ($\vt_t + \vu \cdot \Grad \vt =0$, $\vu = \Curl \Lambda^{-\beta}\vt$ for  $\beta \in [1,2]$) 
possesses such a sufficiently regular solution.


\subsection{Fractional calculus} \label{sec:4}
The fractional Laplace operator
$\Lambda^\alpha$ occurring in \eqref{eq:qg}
is defined using Fourier transform, namely
\begin{equation*}
	\Lambda^\alpha f = \mathcal{F}^{-1}\left(|\dott|^\alpha \mathcal{F}(f)\right).
\end{equation*}
Our normalization of the Fourier transform reads
\begin{equation*}
	 \mathcal{F}f(z)=\int_{\R^2}
         f(\xi)e^{-2\pi i z\cdot\xi}d\xi, \quad  
	 \mathcal{F}^{-1}f(\xi)=\int_{\R^2}
         f(z)e^{2\pi iz\cdot\xi}dz.
\end{equation*}
In the upcoming analysis, we will need a
different representation of $\Lambda^\alpha$. 
Since we require  $\alpha \in (0,2)$,  \cite{Cordoba} provides the
identity 
\begin{equation*}
	\Lambda^\alpha f =C_\alpha\, \pvint_{\R^2}\frac{\delta_z
          f}{|z|^{2+\alpha}}~dz, \quad
        C_\alpha=\frac{\Gamma(1-\alpha/2)}{\pi 2^\alpha \Gamma(\alpha/2)},
\end{equation*}
for all $f$ in the Schwartz class (in particular all $f \in C_0^\infty$). 
Here, $ \pvint \cdots dz$ denotes the principal 
value integral, and we have introduced the notation
\begin{equation*}
	\delta_z f(\xi) = f(\xi+z) - f(\xi), 
	\qquad \xi,z\in\R^2.
\end{equation*}
We will make use of the following 
Leibniz-like formula ``with remainder'':
\begin{equation*}
	\begin{split}
		G^\alpha(f,g) :=\Lambda^\alpha (f g ) 
		- f \Lambda^\alpha g - g \Lambda^\alpha f.
	\end{split}
\end{equation*}	

By adding and subtracting, we see that
\begin{equation}\label{eq:addsubtract}
	\begin{split}
		\delta_z(fg) = g\, \delta_z f + f \,\delta_zg + \delta_z f \,\delta_z g.
	\end{split}
\end{equation}
Multiplying \eqref{eq:addsubtract} with $|z|^{-2-\beta}$ and integrating over $z$ 
provides the following representation of $G^\beta$.	
\begin{lemma}[Leibniz formula]\label{id:product}
If $f$ and $g$ are sufficiently smooth functions 
the following identity holds pointwise, 
\begin{equation*}
  G^\alpha(f,g)=	\Lambda^\alpha (f g )- f \Lambda^\alpha g - g
  \Lambda^\alpha f =
\begin{cases}C_\alpha\, \pvint_\R\frac{\delta_z
      f\, \delta_z g}{|z|^{2+\alpha}}~dz, & \text{for $\alpha \in
      (1,2)$}, \\
2 \nabla f\cdot\nabla g, &  \text{for $\alpha= 2$.}
\end{cases}
\end{equation*}
\end{lemma}

In our analysis, we will make several applications of the above Leibniz rule.
The following proposition provides an $L^p$ estimate for the term $G^\alpha$.
It is a variation of a result due to Constantin \cite{Constantin}.
\begin{proposition}\label{lem:constantin2}
Let $1 < p < \infty$.  For each fixed $\alpha \in (0, 2]$, there is 
a constant $C$ depending on $\alpha$ such that
\begin{equation*}
	\begin{split}
		&\|G^\alpha(f,g)\|_{L^p} \leq 
		C\left(\|\Grad f\|_{L^{\infty}}\|\Grad g\|_{L^{p}}\right)^{\alpha/2}
		\left(\|f\|_{L^{\infty}}\|g\|_{L^{p}}\right)^{1-\alpha/2},			
	\end{split}
\end{equation*}
for all $f$, $g$ $\in C_0^\infty$.
\end{proposition}
\begin{proof}
The result follows directly from the H\"older inequality when $\alpha=2$. 
We may thus assume that $\alpha \in (0, 2)$. 
Let us write $G^\alpha(f,g)$ as the sum of two parts:
\begin{equation*}
	\begin{split}
	G^\alpha(f,g) &= \pvint_{|z|\leq r} 
	\frac{\delta_z f \delta_z g}{|z|^{2+\alpha}}~dz
	+ \int_{|z|>r }\frac{\delta_z f \delta_z g}{|z|^{2+\alpha}}~dz \\
	&:= J^\alpha_r(f,g) + K^\alpha_r(f,g).
	\end{split}
\end{equation*}
(i)  We commence by estimating the first term:
\begin{equation*}
	\begin{split}
	\|J^\alpha_r(f,g)\|_{L^p}^p 
	&= \int_{\R^2} \left|\int_{|z|\leq r}\frac{\delta_z f \delta_z g}{|z|^{2
	+\alpha}}~dz\right|^p~dx \\
	& \leq \int_{\R^2} \left(\int_{|z|\leq r}
	\frac{|\delta_z f \delta_z g|}{|z|^{2+\alpha}}~dz\right)^p~dx \\
	&= \int_{\R^2} \int_{|q_1|\leq r}\ldots \int_{|q_p|\leq r} \prod_{i=1}^p 
	\frac{|\delta_{q_i}f \delta_{q_i}g|}{|q_i|^{2+\alpha}}~dq_1\ldots dq_p\, dx. 
	\end{split}
\end{equation*}
Since $\sum_{i=1}^p \frac{1}{q_i} = 1$, we can apply 
the generalized H\"older inequality to obtain
\begin{equation*}
\begin{split}
\|J^\alpha_r(f,g)\|_{L^p}^p
&\leq \int_{|q_1|\leq r}\ldots \int_{|q_p|\leq r}
\prod_{i=1}^p |q_i|^{-2-\alpha}\|\delta_{q_i}f\delta_{q_i}g\|_{L^{p}}~dq_1\ldots dq_p \\
&\leq \int_{|q_1|\leq r}\ldots \int_{|q_p|\leq r}
\prod_{i=1}^p |q_i|^{-2-\alpha}\|\delta_{q_i}f\|_{L^{\infty}}\|\delta_{q_i}g\|_{L^{p}}~dq_1\ldots dq_p \\	
&\leq \int_{|q_1|\leq r}\ldots \int_{|q_p|\leq r}
\prod_{i=1}^p |q_i|^{-\alpha}\|\Grad f\|_{L^{\infty}}\|\Grad g\|_{L^{p}}~dq_1\ldots dq_p \\
&= \|\Grad f\|^p_{L^{\infty}}\|\Grad g\|^p_{L^{p}}\left(\int_{|z|\leq r}|z|^{-\alpha}~dz \right)^p \\
&\leq  r^{p(2-\alpha)} C\|\Grad f\|^p_{L^{\infty}}\|\Grad g\|^p_{L^{p}},
\end{split}
\end{equation*}
where we have used that 
$\|\delta_z f\|_{L^p} \leq |z|\|\Grad f\|_{L^{p}}$, for all $f \in W^{1,p}$. \\
(ii) A calculation similar to the previous yields
\begin{equation*}
	\begin{split}
		\|K^\alpha_r(f,g)\|_{L^p}^p & 
		\leq \int_{|q_1|> r}\ldots \int_{|q_p|> r}
		\prod_{i=1}^p |q_i|^{-2-\alpha}
		\|\delta_{q_i}f\|_{L^{\infty}}\|\delta_{q_i}g\|_{L^{p}}~dq_1\ldots dq_p \\
		&\leq 4^p\int_{|q_1|> r}\ldots \int_{|q_p|> r}
		\prod_{i=1}^p |q_i|^{-2-\alpha}\| f\|_{L^{\infty}}\| g\|_{L^{p}}~dq_1\ldots dq_p \\
		&\leq C\| f\|_{L^{\infty}}^p\| g\|_{L^{p}}^p 
		\left(\int_{|z|>r}|z|^{-2-\alpha}~dz \right)^p \\
		&\leq r^{-p\alpha}C\| f\|_{L^{\infty}}^p\| g\|_{L^{p}}^p.
	\end{split}
\end{equation*}
Optimizing in $r$ yields the result.
\end{proof}

In order to work efficiently in our Sobolev spaces we will need to
provide some standard definitions, mostly to fix the notation. Let $l$
denote a two-dimensional multi-index, i.e., $l=(l_1,l_2)$, $l_j\in\N_0$. Then we write
\begin{equation*}
  D^l f=\Grad^{l} f =\frac{\partial^{\abs{l}}f}{\partial
    x^{l_1}\partial y^{l_2}}, \qquad \abs{l}=l_1+l_2.
\end{equation*}
If $\ell\in\N$ we let
\begin{equation*}
 \Grad^{\ell} f = \{\Grad^{l} f\mid |l| = \ell\}
\end{equation*}
and 
\begin{equation*}
 \Grad^{\ell} f:\Grad^{\ell} g = \sum_{\substack{l \\ \abs{l}\le \ell}}  \Grad^{l} f\, \Grad^{l} g.
\end{equation*}
We will be working the Sobolev spaces
\begin{equation*}
  H^k=H^k(\R)=\{f\in\mathcal S' \mid (1+\abs{\xi}^2)^{k/2}\mathcal F(
  f(\xi))\in L^2(\R)\}, \quad k\in\R
\end{equation*}
(where $\mathcal S'$ denotes the set of tempered distributions). If $k$ is a natural number, 
$H^k$ is the standard Sobolev space with inner product and norm given by
\begin{equation*}
  \langle f,g\rangle_{H^k} 
=\sum_{\ell = 0}^k\langle \Grad^\ell f,\Grad^\ell g\rangle_{L^2},
\quad \norm{f}_{H^k}=  \langle f,f\rangle_{H^k}^{1/2},
\end{equation*}
where we have introduced
\begin{equation*}
 \langle \Grad^\ell f,\Grad^\ell g\rangle_{L^2}= 
 \sum_{\abs{l} = \ell} \langle D^l f,D^l g\rangle_{L^2}.
\end{equation*} 

In our analysis, we will apply the following corollary of the previous proposition.

\begin{corollary}\label{cor:com}
For each fixed $\alpha \in(0, 2]$ and  integer $k \geq 3$, there is a constant $C$ depending 
on $\alpha$ such that
\begin{equation*}
	\|G^\alpha(f,g)\|_{H^k} \leq C\|f\|_{H^{k+1}}\|g\|_{H^{k+1}},
\end{equation*}
for all $f$, $g$ $\in H^{k+1}$.
\end{corollary}

\begin{proof}
Since $C_0^\infty$ is dense in $H^{k+1}$, we may
assume that $f$, $g$ $\in C_0^\infty$.
By definition,
\begin{equation*}
	\|G^\alpha(f,g)\|^2_{H^k} = \sum_{s=0}^k \|\Grad^s G^\alpha(f,g)\|_{L^2}^2.
\end{equation*}
Our strategy is to prove the desired estimate for each of terms in the 
sum separately. For this purpose, we let $s=0,\ldots, k$ be arbitrary 
and consider an arbitrary component of $\Grad^s G^\alpha(f,g)$. Using 
Lemma \ref{id:product} and the standard Leibniz rule, we deduce
\begin{equation}\label{eq:eqeq}
	\begin{split}
		&\frac{\partial^s}{\partial x^l \partial y^{s-l}}G^\alpha(f,g) \\
		&\qquad 
		= \text{p.v.}\int_{\R^2}\frac{1}{|z|^{2+\alpha}}\sum_{n=0}^l
		\sum_{m=0}^{s-l} \left( l \atop n\right)\left(s-l \atop m\right) \\
		&\qquad \qquad\qquad \qquad \times 
		  \left(\frac{\partial^{m+n}}{\partial x^n\partial y^m}\partial_z f\right)
		  \left(\frac{\partial^{s-m-n}}{\partial x^{l-n}\partial y^{s-l-m}}\partial_z g\right)~dz \\
		&\qquad
		= \sum_{n=0}^l\sum_{m=0}^{s-l} \left( l \atop n\right)\left(s-l \atop m\right)
		G^\alpha\left(\frac{\partial^{m+n}}{\partial x^n\partial y^m}f, 
		\frac{\partial^{s-m-n}}{\partial x^{l-n}\partial y^{s-l-m}} g \right).
	\end{split}
\end{equation}
For $m$, $n$ such that $2 \leq m+n \leq k$, Proposition \ref{lem:constantin2} can 
applied to conclude
\begin{equation*}
	\begin{split}
		&\left\|G^\alpha\left(\frac{\partial^{m+n}}{\partial x^n\partial y^m}f, 
		\frac{\partial^{s-m-n}}{\partial x^{l-n}\partial y^{s-l-m}} g \right)\right\|_{L^2} \\
		&\qquad \qquad 
		\leq C\|f\|_{H^{k+1}}\|\Grad^{s-m-n+1}g\|_{L^\infty}^{\alpha/2}\|\Grad^{s-m-n}g\|_{L^\infty}^{1-\alpha/2}
		\leq C\|f\|_{H^{k+1}}\|g\|_{H^{k+1}}.
	\end{split}
\end{equation*}
Conversely, if $m$,$n$ is such that $0 \leq m+n \leq 1$, 
Proposition \ref{lem:constantin2} allow us deduce
\begin{equation*}
	\begin{split}
		&\left\|G^\alpha\left(\frac{\partial^{m+n}}{\partial x^n\partial y^m}f, 
		\frac{\partial^{s-m-n}}{\partial x^{l-n}\partial y^{s-l-m}} g \right)\right\|_{L^2} \\
		&\qquad \qquad 
		\leq C\|g\|_{H^{k+1}}\|\Grad^{m+n+1}f\|_{L^\infty}^{\alpha/2}
		\|\Grad^{m+n}f\|_{L^\infty}^{1-\alpha/2} \\
		&\qquad \qquad
		\leq C\|g\|_{H^{k+1}}\|f\|_{H^{4}} \leq C\|g\|_{H^{k+1}}\|f\|_{H^{k+1}},
	\end{split}
\end{equation*}
since $k \geq 3$. Now, by taking the $L^2$ norm on both sides of \eqref{eq:eqeq} and applying the 
previous calculations, we gather
 \begin{equation*}
 	\begin{split}
 		&\left\|\frac{\partial^s}{\partial x^l \partial y^{s-l}}G^\alpha(f,g)\right\|_{L^2}  
		\leq C\|f\|_{H^{k+1}}\|g\|_{H^{k+1}}.
	\end{split}
\end{equation*}

Hence, any given component of $\Grad^s G^\alpha(f,g)$ satisfies the desired bound. Thus,
\begin{equation*}
	\|G^\alpha(f,g)\|^2_{H^k} = \sum_{s=0}^k \|\Grad^s f\|_{L^2}^2 
	\leq C\|f\|_{H^{k+1}}\|g\|_{H^{k+1}},
\end{equation*}
which completes our proof.
\end{proof}

\subsection{Two auxiliary lemmas}
We will make heavily use of the following two 
lemmas throughout the paper.
Their proofs are technical and tedious, but straightforward.
For this reason, proofs are deferred to the appendix. 
\begin{lemma}\label{lem:nonlineardiv}
Let $k\geq 6$ be an integer. Then, 
\begin{equation}
	\sum_{s=0}^k\left|\int_{\R^N} \Grad^s(\Grad f\cdot \Curl \Lambda^{-\beta}f):\Grad^s f~dx\right| \leq C\|f\|_{H^{k-2}}\|f\|_{H^k}^2,
\end{equation}	
for all $f \in H^k$.
\end{lemma}
\begin{lemma}\label{lem:lineardiv}
Let $k\geq 4$ be an integer. The following estimates hold  
	\begin{align}\label{eq:ldiv1}
		\sum_{s=0}^k\left|\int_{\R^N} \Grad^s \left(\Grad f \cdot \Curl \Delta^{-\beta} g  \right):\Grad^s f~dx\right| &\leq C\|g\|_{H^k}\|f\|_{H^k}^2,  \quad  f,g\in H^k,\\
\label{eq:ldiv2}
			\sum_{s=0}^k\left|\int_{\R^N}\Grad^s\left(\Grad g \cdot \Curl \Lambda^{-\beta}f\right):\Grad^s f~dx\right| &\leq C\|g\|_{H^{k+1}}\|f\|_{H^k}^2, \quad  f\in H^k, \,  g\in H^{k+1}.
	\end{align}
\end{lemma}

\section{Godunov splitting}\label{sec:godunov}
In this section we prove rigorously the expected linear 
rate of convergence for Godunov splitting. 
As part of the proof we also show that this splitting 
method is well-defined and that it 
produces regular approximations. 
The results are valid under a condition 
on the length of the time step $\Delta t$.

We begin by precisely defining Godunov splitting for \eqref{eq:qg}. For this purpose, 
write \eqref{eq:qg} as
\begin{equation*}
	\vt_t = C(\vt), \qquad C(\vt) = A(\vt) + B(\vt), 
\end{equation*}
\begin{equation*}
	B(\vt) = -\Curl \Lambda^{-\beta}\vt \cdot \Grad \vt, \qquad A(\vt) 
	= - \Lambda^{\alpha}\vt.
\end{equation*}
Using the operators $A$ and $B$, we define the solution 
operators $\Phi_A$ and $\Phi_B$ as the solutions to
\begin{equation*}
	\begin{split}
		\partial_t\Phi_A(t,\vt_0) &= A(\Phi_A(t,\vt_0)), \quad \Phi_A(0,\vt_0)= \vt_0, \\
		\partial_t\Phi_B(t,\vt_0) &= B(\Phi_B(t,\vt_0)), \quad \Phi_B(0,\vt_0)= \vt_0.
	\end{split}
\end{equation*}

The Godunov splitting method is classically defined as follows:
For $\Delta t> 0$ given,  construct a sequence 
$\{\vt^n, \vt^{n+1/2}\}_{n=0}^{\lfloor T/\Delta t\rfloor}$ of approximate solutions 
to \eqref{eq:qg} by the following procedure:
Let $\vt^0 = \vt_0$ and determine inductively
\begin{equation*}
\vt^{n+1/2} = \Phi_B(\Delta t,\vt^{n}),\quad    
\vt^{n+1} = \Phi_A(\Delta t, \vt^{n+1/2}), 
\quad n=0, \ldots, \lfloor T/\Delta t\rfloor-1.    
\end{equation*}

To facilitate the convergence analysis, we will need 
a different definition of the Godunov method. Our definition 
can be seen as an extension of the splitting 
solution $\{\vt^n, \vt^{n+1/2}\}_n$ to all of $[0,T]$.
The most used method of extension is to let ``time run twice as fast"
in each of the sub-intervals $[t_n, t_{n+1/2}]$
and $[t_{n+1/2}, t_{n+1}]$, where as usual 
$t_r=r\Dt$ for $r\in [0,\infty)$, thus obtaining
\begin{equation*}
	\vt_{\Delta t}(t) = 
	\begin{cases}
		\Phi_B(2(t-t_n), \vt^n), & t \in [t_n, t_{n+1/2}), \\
		\Phi_A(2(t-t_{n+1/2}), \vt^{n+1/2}), & t \in [t_{n+1/2}, t_{n+1}).
	\end{cases}
\end{equation*}
Although it appears to be a natural extension, it 
does not seem to be appropriate for our purpose. 
Instead, we will follow the approach taken in the paper \cite{Holden:tao}, 
and introduce two time variables instead of one. 
\begin{figure}[tbp]
  \centering
  \includegraphics[width=0.5\linewidth]{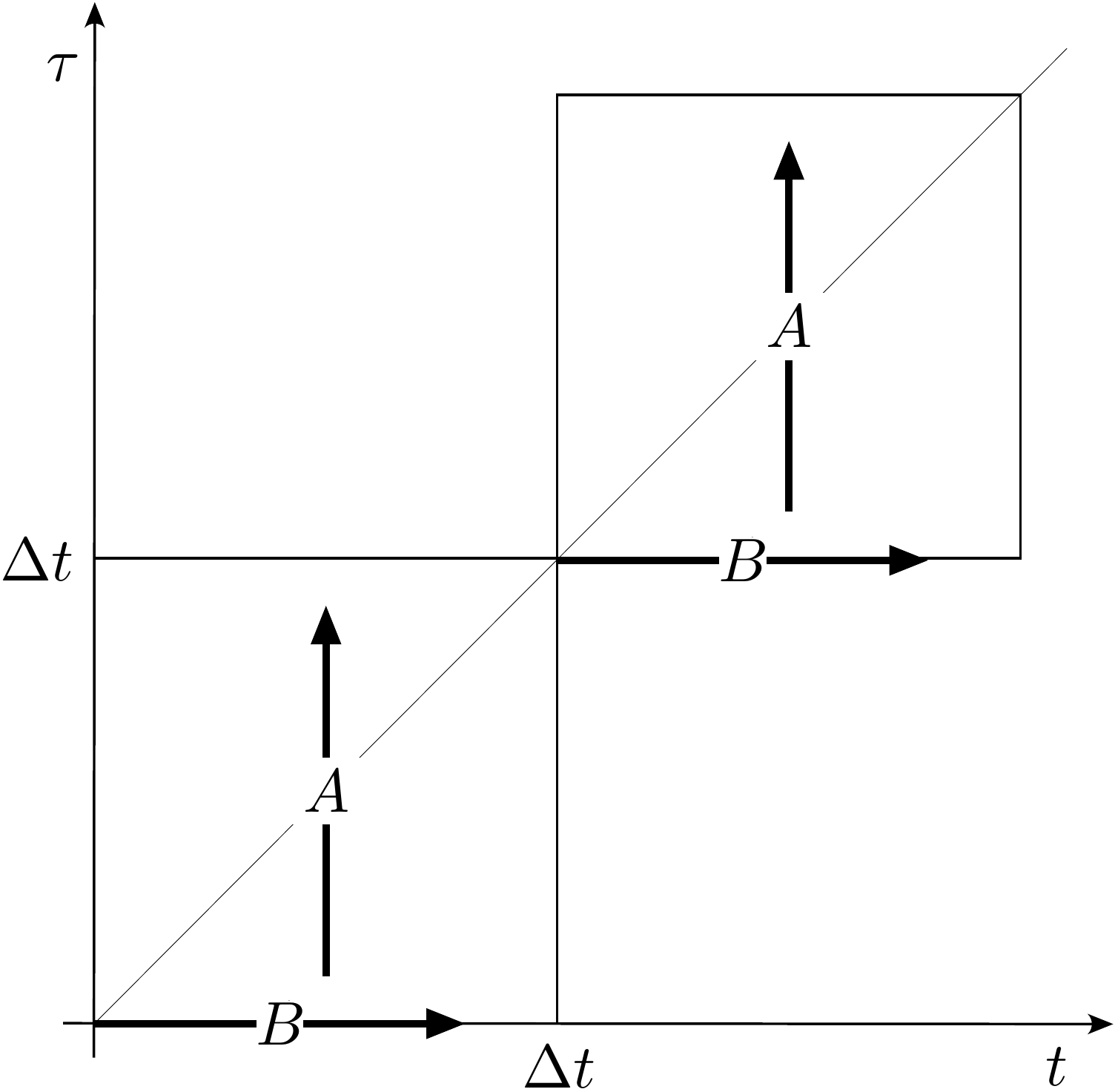}
  \caption{A schematic view of Godunov splitting.}
  \label{fig:GodunovSplit}
\end{figure}

Our Godunov splitting method is given by the following definition.

\begin{definition}[Godunov splitting]\label{def:godunov}
With $\Delta t>0$ given, define the domain
\begin{equation*}
	\Om_{\Delta t} 
	= \bigcup_{n=0}^{\lfloor T/\Delta t\rfloor-1}[t_n, t_{n+1}]
	\times [t_n, t_{n+1}].
\end{equation*}
The time-continuus Godunov splitting solution 
$\vvt:\Om_{\Delta t} \to \R$ is defined as 
the solution to 
\begin{equation}\label{eq:godunov-2}
	\begin{split}
		&\vvt(0,0) = \vt_0, \\
		& 
		\begin{cases}
			\vvt_t(t,t_n) = B(\vvt(t,t_n)), & t \in (t_n, t_{n+1}], \\
			\vvt|_{t=t_n}=\theta^n;
		\end{cases}\\
		& 
		\begin{cases}
			\vvt_\tau(t, \tau) =  A(\vvt(t, \tau)),& (t,\tau) 
			\in [t_n, t_{n+1}]\times (t_n,t_{n+1}],\\
			\vvt|_{\tau=t_n}=\vvt(t,t_n), & t\in  [t_n, t_{n+1}],
		 \end{cases}
	\end{split}
\end{equation}
cf.~Figure \ref{fig:GodunovSplit}.
\end{definition}

Observe that
\begin{equation*}
 \vvt(t_n,t_n)=\theta^n, \quad  n=0, \ldots, \lfloor T/\Delta t\rfloor.
\end{equation*}
Thus, $\vvt(t,t)$ is an extension of $\{\vt^n\}_n$ to all of $[0,T]$.

To measure the error, we will use the function
\begin{equation*}
	e(t) =   \vvt(t,t)-\vt(t),
\end{equation*}
where $\vt$ is the (smooth) solution of \eqref{eq:qg}.

It is not trivial trivial to obtain the existence of a splitting solution
$\vvt$ in the sense of Definition \ref{def:godunov}.
Since the best available existence result 
we have for the hyperbolic step is local-in-time, 
it is unclear if we can iterate the steps up to any given large time $T$.
In fact, well-posedness of the method is one of our main results.

The following theorem is our main result in this section.
\begin{theorem}\label{thm:main}
Suppose $\theta_0 \in H^k$ with $6 \leq k \in \mathbb{N}$, and  that $\alpha$, $\beta$ $\in [1,2]$. 
Then, for $\Delta t$ sufficiently small, we have the following:
\begin{enumerate}
	\item The Godunov splitting approximation 
	$\vvt(t,\tau)$ is well-defined. Moreover, $\vvt(t,\tau)$ belongs to $C([0,T]^2;H^k)$. 
	\item The error $e(t):= \|\vvt(t,t) - \vt(t)\|_{H^{k-2}} $ satisfies
	$$
	\|e(t)\|_{H^{k-2}} \leq tC\Delta t\|\theta_0\|_{H^k}^2.
	$$
\end{enumerate}
\end{theorem}

Theorem \ref{thm:main} will be an outcome of 
the results proved in the subsections below 
(Subsection \ref{sec:prf} will bring the pieces together).

\subsection{Evolution equations for the error}
To prove Theorem \ref{thm:main}, we will analyze 
a set of evolution equations  satisfied by the error $e$, see  \cite{Holden:tao}, which 
we now derive.

We shall need the following Taylor expansion satisfied 
by an operator $E$:
\begin{align*}
	E(f+g) &= E(f) + dE(f)[g] 
	+ \int_0^1 (1-\gamma)d^2E(f+\gamma g)[g]^2~d\gamma.
\end{align*}
Using the definition of $\vvt$ and the above Taylor formula, we deduce
\begin{equation}\label{eq:err1}
	\begin{split}
		e_t - dC(\vt)[e] 
		&= \vvt_t + \vvt_\tau - \vt_t - dA(\vt)[e] - dB(\vt)[e] \\
		&= \vvt_t + A(\vvt) - (A+B)(\vt) - dA(\vt)[e] - dB(\vt)[e] \\
		&= \vvt_t - B(\vvt) + (A(\vvt) - A(\vt) - dA(\vt)[e]) \\
		&\qquad \qquad +  (B(\vvt) - B(\vt) - dB(\vt)[e]) \\
		&= F(t,t) + \int_0^1 (1-\gamma)d^2C(\vt+ \gamma e)[e]^2~d\gamma,
	\end{split}
\end{equation}
where we have introduced the ``forcing" term
\begin{equation*}
	F(t, \tau) = \vvt_t(t, \tau)  - B(\vvt(t,\tau)).
\end{equation*}
By direct calculation,
\begin{equation}\label{eq:err3}
	\begin{split}
		F_\tau - dA(\vvt)[F] &= v_{t\tau} - B(\vvt)_\tau - dA(\vvt)[\vvt_t - B(\vvt)] \\
		&= A(\vvt)_t - dB(\vvt)[\vvt_\tau] - dA(\vvt)[\vvt_t] + dA(\vvt)[B(\vvt)] \\
		&= dA(\vvt)[\vvt_t] - dB(\vvt)[A(\vvt)] - dA(\vvt)[\vvt_t] + dA(\vvt)[B(\vvt)] \\
		&= [A,B](\vvt),
	\end{split}
\end{equation}
where we have defined the commutator
\begin{equation*}
	\begin{split}
		[A,B](f) = dA(f)[B(f)] - dB(f)[A(f)].
	\end{split}
\end{equation*}

For the fluid equation \eqref{eq:qg}, we have that
\begin{equation*}
	\begin{split}
		A(f) &= -\Lambda^\alpha f, \\
		dA(f)[g] &= -\Lambda^\alpha g, \\
		d^2A(f)[g,h] &= 0,
	\end{split}
\end{equation*}
and
\begin{equation*}
	\begin{split}
		B(f) &= -\Grad f \cdot \Curl \Lambda^{-\beta}f, \\
		dB(f)[g] &= -\Grad f\cdot \Curl \Lambda^{-\beta}g 
		- \Grad g\cdot \Curl \Lambda^{-\beta}f, \\
		d^2B(f)[g,h] &= -\Grad h\cdot \Curl \Lambda^{-\beta}g 
		- \Grad g\cdot \Curl \Lambda^{-\beta}h.
	\end{split}
\end{equation*}
The Leibniz formula (Lemma \ref{id:product}) yields
\begin{equation*}
	\begin{split}
		[A,B](f) &= 
		\Lambda^\alpha\left(\Grad f \cdot \Curl \Lambda^{-\beta}f\right) \\
		&\qquad - \Grad f\cdot \Curl \Lambda^{-\beta}\Lambda^\alpha f
		- \left(\Grad \Lambda^\alpha f\right)\cdot \Curl \Lambda^{-\beta}f \\
		&= G^\alpha(\Grad f, \Curl \Lambda^{-\beta}f).
	\end{split}
\end{equation*}
Hence, for the fluid equation \eqref{eq:qg}, equations 
\eqref{eq:err1}--\eqref{eq:err3} read:
\begin{align}
F &= \vvt_t + \Grad \vvt \cdot \Curl \Lambda^{-\beta}\vvt,  \notag 
\\ F_\tau + \Lambda^\alpha F&= G^\alpha(\Grad \vvt, \Curl \Lambda^{-\beta}\vvt), 
\label{eq:feq}\\
e_t + \Lambda^\alpha e &= F -\Grad e\cdot \Curl \Lambda^{-\beta}\theta 
-\Grad \theta\cdot \Curl \Lambda^{-\beta}e
-\Grad e\cdot \Curl \Lambda^{-\beta}e. \label{eq:weq}
\end{align}
The equations \eqref{eq:feq} and \eqref{eq:weq} constitute
our main tool for proving Theorem \ref{thm:main}.

\subsection{Estimates on the error}

Lemma \ref{lem:w} below provides an estimate on the error $e$ 
that will be a key ingredient in the proof of Theorem \ref{thm:main}.
The estimate depends on two auxiliary results (Lemmas \ref{lem:higher} and \ref{lem:F}) that 
we  first prove.

The first auxiliary results gives perhaps 
the most fundamental property of our splitting solution.
It states that if the splitting solution is in
$H^k$, $k\geq 6$, then 
it is actually in $H^{k+2}$. 

For notational convenience, let $\Om_{\Delta t}^{t,\tau}$ 
denote the set of all times prior to $(t, \tau)$:
\begin{equation*}
\Om_{\Delta t}^{t,\tau} = \{(s,\sigma) \in \Om_{\Delta t} \mid 0 \leq
s \leq t, ~0 \leq \sigma \leq \tau\}
=\Om_{\Delta t}\cap\big([0,t]\times [0,\tau]\big) .
\end{equation*}

\begin{lemma}\label{lem:higher}
Assume that $\alpha$, $\beta$ $\in [1,2]$, 
 $6 \leq k \in \mathbb{N}$, and let $\vvt$ be the Godunov 
approximation to \eqref{eq:qg} in the sense of Definition \ref{def:godunov} and
\eqref{eq:godunov-2}.

If for some $(s,\sigma) \in \Om_{\Delta}$,
\begin{equation}\label{eq:gamma}
	\|\vvt(s,\sigma)\|_{H^{k-2}}\leq \gamma,\quad (s,\sigma) \in \Om_{\Delta t}^{t,\tau},
\end{equation}
then
\begin{equation*}
	\|\vvt(s,\sigma)\|_{H^{k}} \leq e^{C\gamma s}\|\vt_0\|_{H^k}, 
\quad (s,\sigma) \in \Om_{\Delta t}^{t,\tau}. 
\end{equation*}
\end{lemma}
\begin{proof}
By direct calculation we see that
\begin{equation*}
	\begin{split}
		\partial_\tau\frac12 \|\vvt (t,\tau)\|_{H^k}^2 
		&= \sum_{\ell=0}^k\int_{\R^2} \Grad^\ell \vvt_\tau: \Grad^\ell \vvt~dx \\
		&= - \sum_{\ell=0}^k\int_{\R^2} \left|\Grad^\ell\Lambda^{\alpha/2}\vvt\right|^2~dx 
		= - \|\vvt\|_{H^{k + \alpha/2}}^2.
	\end{split}
\end{equation*}
Hence, for $(t, \tau) \in [t_n,t_{n+1}] 
\times (t_n, t_{n+1}]$, $n=0,\ldots, \lfloor T/\Dt\rfloor-1$,
$$
\|\vvt (t,\tau)\|_{H^k}^2  \leq \|\vvt (t,t_n)\|_{H^k}^2.
$$

We now calculate a bound on the $H^{k}$ norm of $\vvt(t,t_n)$. By definition,
\begin{equation}\label{eq:startburger}
	\begin{split}
		\partial_t\frac12 \|\vvt (t,t_n)\|_{H^{k}}^2 
		&= \sum_{\ell=0}^k\int_{\R^2} \Grad^\ell \vvt_t : \Grad^\ell\vvt~dx \\ 
		&= -\sum_{\ell=0}^k\int_{\R^2}  
		\Grad^\ell\left(\Grad \vvt\cdot \Curl \Lambda^{-\beta}\vvt\right):\Grad^\ell\vvt~dx.
	\end{split}
\end{equation}
Applying Lemma \ref{lem:nonlineardiv} to this identity, we see that
\begin{equation}
	\partial_t\frac12 \|\vvt (t,t_n)\|_{H^{k}}^2 \leq C\|\vvt(t,t_n)\|_{H^{k-2}}\|\vvt(t,t_n)\|_{H^k}^2.
\end{equation}
Clearly, this allows us to conclude that
\begin{equation*}
	\partial_t \|\vvt(t,t_n)\|_{H^k} \leq C\gamma\|\vvt(t,t_n)\|_{H^k},
\end{equation*}
and integration in time gives 
\begin{equation*}
	\|\vvt(t,t_n)\|_{H^k} \leq e^{C\gamma (t-t_n)}\|\vvt(t_n,t_n)\|_{H^k}.
\end{equation*}
Consequently,
\begin{equation*}
	\|\vvt(t,t_n)\|_{H^k} \leq e^{C\gamma t}\|\vt_0\|_{H^k},
\end{equation*}
which completes the proof. 
\end{proof}

The following lemma is our second auxiliary result.
\begin{lemma}\label{lem:F}
	Let $k$, $t$, and $\tau$ be as in the previous lemma. Then,
	\begin{equation*}
		\|F(s,\sigma)\|_{H^{k-2}} \leq \Delta t\, C_1e^{C_2\gamma s}\|\vt_0\|^2_{H^k}, 
		\quad  (s,\sigma) \in \Om^{t,\tau}_{\Delta t}.
	\end{equation*}
\end{lemma}

\begin{proof}
Applying $\Grad^{s}$ to \eqref{eq:feq}, multiplying componentwise 
with $\Grad^{s}F$, summing over $s=0, \ldots, k-2$, and 
integrating over the domain gives
\begin{equation}\label{eq:F1}
	\begin{split}
		&\partial_\tau \frac{1}{2}\|F\|_{H^{k-2}}^2 + \|F\|_{H^{k-2+\alpha/2}}^2 \\
		&\qquad \qquad
		= \sum_{s=0}^{k-2}\int_{\R^2}\Grad^{s}F: 
		\Grad^{s}G^\alpha(\Grad \vvt, \Curl \Lambda^{-\beta}\vvt)~dx \\
		&\qquad \qquad
		\leq C\|F\|_{H^{k-2}}\|G^\alpha(\Grad \vvt, \Curl
                \Lambda^{-\beta}  \vvt)\|_{H^{k-2}}.
	\end{split}
\end{equation}
An application of Corollary \ref{cor:com} to \eqref{eq:F1} gives 
(here we do not make use of the regularizing effect of $\Lambda^{-\beta}$)
\begin{equation*}
	\begin{split}
		\partial_\tau \frac{1}{2}\|F\|_{H^{k-2}}^2
		\leq C\|F\|_{H^{k-2}}\|\Grad \vvt\|_{H^{k-1}}
		\|\Curl \Lambda^{-\beta}\vvt\|_{H^{k-1}}
		\leq C\|F\|_{H^{k-2}}\|\vvt\|_{H^k}^2.
	\end{split}
\end{equation*}
In view of Lemma \ref{lem:higher}, this means that
\begin{equation*}
	\partial_\tau \|F\|_{H^{k-2}} \leq C_1 e^{C_2\gamma t}\|\vt_0\|_{H^k}^2.
\end{equation*}
Since $F(t,t_n) = 0$ integration in time gives
\begin{equation*}
	\|F\|_{H^{k-2}} \leq \Delta t\, C_1 e^{C_2\gamma t}\|\vt_0\|_{H^k}^2,
\end{equation*}
which concludes the proof.
\end{proof}

The next lemma will be the key ingredient in the proof of Theorem \ref{thm:main}.
\begin{lemma}\label{lem:w}
 Let $k$, $t$, and $\tau$ be as in Lemma \ref{lem:higher}. Then,
\begin{equation*}
	\|e(s)\|_{H^{k-2}} \leq s\Delta t\, C_1e^{c\gamma t}\|\vt_0\|_{H^k}^2, \quad 0 < s \leq t.
\end{equation*}
\end{lemma}

\begin{proof}
Applying $\Grad^s$ to \eqref{eq:weq}, multiplying componentwise 
with $\Grad^s e$, summing over $s=0, \ldots, k-2$, and 
integrating gives
\begin{equation}\label{eq:w1}
	\begin{split}
		&\frac{1}{2}\partial_t\|e\|_{H^s}^2 + \|e\|_{H^{s+\alpha/2}}^2 \\
		&\qquad = \sum_{s=0}^{k-2}\int_{\R^2}\Big[ \Grad^{s} F\,: \Grad^{s} e 
		+\Grad^s\left(\Grad e\cdot \Curl \Lambda^{-\beta}e\right):\Grad^s e \\
		&\qquad \qquad - \Grad^s\left(\Grad \vt\cdot \Curl \Lambda^{-\beta}e
		+\Grad e\cdot \Curl \Lambda^{-\beta}\vt \right):\Grad^s e\Big]dx \\
		&\qquad \leq \|e\|_{H^{k-2}}\|F\|_{H^{k-2}}
		+\sum_{s=0}^{k-2}\int_{\R^2}|\, \Grad^s\left(\Grad e\cdot 
		\Curl \Lambda^{-\beta}e\right):\Grad^s e \\
		&\qquad \qquad - \sum_{s=0}^{k-2}
		\Grad^s\left(\Grad \vt\cdot \Curl \Lambda^{-\beta}e
		+\Grad e\cdot \Curl \Lambda^{-\beta}\vt
                      \right):\Grad^s e\, |~dx.
	\end{split}
\end{equation}
An application of Lemma \ref{lem:nonlineardiv} provides the estimate
\begin{equation}\label{eq:w2}
	\begin{split}
		&\sum_{s=0}^{k-2}\int_{\R^2}\Grad^s\left(\Grad e\cdot 
		\Curl \Lambda^{-\beta}e\right):\Grad^s e~dx 
		\leq C \|e\|_{H^{k-2}}^3 \\
		&\qquad\leq C\|e\|_{H^{k-2}}^2\left(\|\vvt\|_{H^{k-2}} + \|\vt\|_{H^{k-2}}\right) 
		 \leq C\|e\|_{H^{k-2}}^2\left(\gamma + \|\theta_0\|_{H^{k-2}}\right),
	\end{split}
\end{equation}
where we have also used that $\|\vvt\|_{H^{k-2}} \leq \gamma$.

By virtue of Lemma \ref{lem:lineardiv}, we also have the estimate
\begin{equation}\label{eq:w3}
	\begin{split}
		&\sum_{s=0}^{k-2}\int_{\R^2} \Grad^s\left(\Grad \vt\cdot \Curl \Lambda^{-\beta}e
		 		+\Grad e\cdot \Curl \Lambda^{-\beta}\vt \right):\Grad^s e~dx \\		
		&\qquad \leq C\left(\|\vt\|_{H^{k-1}}\|e\|_{H^{k-2}}^2 + \|\vt\|_{H^{k-2}}\|e\|_{H^{k-2}}^2\right).
	\end{split}
\end{equation}

Combining \eqref{eq:w3},
 \eqref{eq:w2}, and applying the result together with Lemma \ref{lem:F} to \eqref{eq:w1} enable us to conclude that
\begin{equation*}
	\partial_t\|e\|_{H^{k-2}} 
	\leq \Delta t\, C_1 e^{C_2\gamma t}\|\vt_0\|_{H^k}^2
	+  C(\gamma + \|\theta^0\|_{H^k})\|e\|_{H^{k-2}}.
\end{equation*}
An application of the Gronwall inequality (recalling that $e(0) = 0$) yields the result.
\end{proof}

\subsection{Proof of Theorem \ref{thm:main}}
\label{sec:prf}

We will make use of the following bootstrap lemma (cf.~\cite[Proposition 1.21]{Tao}):

\begin{lemma}\label{lem:bootstrap}
For each $(t,\tau) \in \Om_{\Delta t}$, suppose that we have two statements, a ``hypothesis" 
$\vc{H}(t, \tau)$ and a ``conclusion" $\vc{C}(t, \tau)$. Suppose that we can verify the 
following assertions:
\begin{enumerate}
	\item If $\vc{H}(t, \tau)$ is true for some  $(t,\tau) \in \Om_{\Delta t}$, then 
	$\vc{C}(t, \tau)$ is  also true. 
	\item If $\vc{C}(t, \tau)$ is true for some $(t,\tau) \in \Om_{\Delta t}$, then $\vc{H}(t', \tau')$ is 
	also true for all $(t',\tau')$ in a neighborhood of $(t,\tau)$.
	\item If $(t_1, \tau_1)$, $(t_2, \tau_2)$, $\ldots$ is a 
	sequence in $\Om_{\Delta t}$  converging
	to $(t,\tau) \in \Om_{\Delta t}$, and $\vc{C}(t_n, \tau_n)$ is true 
	for all $n$, then $\vc{C}(t, \tau)$ is true.
	\item $\vc{H}(t, \tau)$ is true for at least one  $(t,\tau) \in \Om_{\Delta t}$.
\end{enumerate}
Then $\vc{C}(t, \tau)$ is true for all $(t,\tau) \in \Om_{\Delta t}$.
\end{lemma}

Let $\vc{H}(t, \tau)$ denote the statement
\begin{equation*}
	\|\vvt(s,\sigma)\|_{H^{k-2}} \leq \gamma, \quad  
	(s,\sigma) \in \Om^{t,\tau}_{\Delta t},
\end{equation*}
and $\vc{C}(t, \tau)$  the statement
\begin{equation*}
	\|\vvt(s,\sigma)\|_{H^{k-2}} \leq \frac{\gamma}{2}, 
	\quad  (s,\sigma) \in \Om^{t,\tau}_{\Delta t},
\end{equation*}
where $\gamma$ is some value which will specified below. Let us for a moment 
assume that assertions (1)--(4) of Lemma \ref{lem:bootstrap} are true for 
 $\vc{H}$ and $\vc{C}$. In this case Lemma \ref{lem:bootstrap} tells us 
that $\vc{C}(t, \tau)$ holds for all $(t, \tau) \in \Om_{\Delta t}$. 
Hence, the $H^{k-2}$ norm of the splitting solution never blows up and 
thus the existence part of Theorem \ref{thm:main} follows readily
(a local-in-time solution can be extended to all of $(0,T)$).

\bigskip
Let us now verify assertions (1)--(4) of Lemma \ref{lem:bootstrap} for our 
choice of $\vc{H}(t, \tau)$ and $\vc{C}(t, \tau)$. First, we observe that 
assertions (2) and (3) clearly hold. For assertion (4) to be true, we 
assume that $\gamma$ satisfies
\begin{equation}\label{eq:asmpt}
\|\vvt(0,0)\|_{H^{k-2}} = \|\vt_0\|_{H^{k-2}} \leq \gamma.
\end{equation}
Then it remains to verify assertion (1).
For this purpose, let us assume that $\vc{H}(t, \tau)$ is true for some $(t, \tau) \in \Om_{\Delta t}$.
Lemma \ref{lem:w} can then be applied to obtain the bound
\begin{equation}\label{eq:vf1}
	\|e(s)\|_{H^{k-2}} \leq s(\Delta t)Ce^{c\gamma t}\|\vt_0\|_{H^k}^2, \quad s \in [0,t].
\end{equation}
We need to compare  $\vvt(s,\sigma)$ with the corresponding value on
the diagonal, viz.\  $\vvt(s,s)$. Two cases need to be considered:\\
(i) If $\sigma>s$, we observe that
\begin{equation*}
	\partial_\sigma \frac{1}{2}\|\vvt (s,\sigma)\|_{H^{k-2}}^2 
	= \sum_{s=0}^{k-2}\int_{\R^2} \Grad^{s} \vvt_\sigma: \Grad^{s} \vvt~dx \\
	= - \sum_{s=0}^{k-2}\int_{\R^2}|
	\Grad^{s}\Lambda^{\alpha/2}\vvt|~dx\le 0, 
\end{equation*}
(using that  $\alpha \leq 2$ to make sure the expression is finite)
which implies that  \\ 
$\|\vvt(s,\sigma)\|_{H^{k-2}}\le \|\vvt(s,s)\|_{H^{k-2}}$. \\
(ii) If $\sigma<s$ we first see that
\begin{equation*}
	\begin{split}
		\left |\partial_\sigma \frac{1}{2}\|\vvt (s,\sigma)\|_{H^{k-2}}^2 \right |
		&=\left | \sum_{s=0}^{k-2}\int_{\R^2} 
		\Grad^{s} \vvt_\sigma: \Grad^{s} \vvt~dx \right| \\
		&= \left | \sum_{s=0}^{k-2}\int_{\R^2}
		\Grad^{s}\Lambda^{\alpha}\vvt : \Grad^{s}\vvt~dx\right | 
		\leq \|\vvt(s,\sigma)\|_{H^{k-2}}\|\vvt(s, \sigma)\|_{H^{k}},
	\end{split}
\end{equation*}
(using $\alpha \leq 2$) which implies that 
\begin{equation*}
 	\left |\partial_\sigma \|\vvt (s,\sigma)\|_{H^{k-2}} \right |\le  \|\vvt(s, \sigma)\|_{H^{k}}. 
\end{equation*}
Using this inequality and applying Lemma \ref{lem:higher} we find
\begin{equation*}
\begin{split}
\|\vvt(s, \sigma)\|_{H^{k-2}} 
&\leq \|\vvt(s, s)\|_{H^{k-2}} +
                \int_\sigma^s|\partial_\sigma \|\vvt
                (s,\tilde\sigma)\| |~d\tilde\sigma\\
	&\leq \|\vvt(s, s)\|_{H^{k-2}} 
	+ \Delta t \sup_{\sigma' \in [\sigma,s]}\|\vvt(s,\sigma')\|_{H^k} \\
	&\leq \|\vvt(s, s)\|_{H^{k-2}}
	+ \Delta t\, e^{C\gamma s}\|\theta_0\|_{H^{k}}.
\end{split}
\end{equation*}
Thus in both cases we conclude that
\begin{equation*}
\|\vvt(s, \sigma)\|_{H^{k-2}} \leq \|\vvt(s, s)\|_{H^{k-2}}+
        \Delta t\, e^{C\gamma s}\|\theta_0\|_{H^{k}},
        \quad |s-\sigma|\le \Dt.
\end{equation*}
By applying the previous inequality, adding and subtracting
        $\theta$, and involving \eqref{eq:vf1}, we estimate
\begin{equation}\label{eq:step31}
\begin{split}
	\|\vvt(s,\sigma)\|_{H^{k-2}} 
	&\leq \|e(s)\|_{H^{k-2}} +
                        \|\vt(s)\|_{H^{k-2}}+ \Delta t\, e^{C\gamma s}\|\theta_0\|_{H^{k}} \\
	&\leq C\Delta t\, e^{C\gamma
                          T}(s\|\theta_0\|^2_{H^{k}}+\|\theta_0\|_{H^{k}}\big) + C\|\theta_0\|_{H^{k}} \\
	&\leq \Delta t\, C_1(\gamma)  + C_2.
\end{split}
\end{equation}
Now, we fix $\gamma$ and $\Delta t$ according to
\begin{equation*}
\gamma = 4C_2, \qquad \Delta t \leq \frac{C_2}{C_1(\gamma)},
\end{equation*}
and note that this is not conflict with \eqref{eq:asmpt}.

Then, \eqref{eq:step31} gives
\begin{equation*}
\begin{split}
	\|\vvt(s,\sigma)\|_{H^{k-2}}  \leq 2C_2 = 
	\frac{\gamma}{2},\quad (s,\sigma) \in \Om^{t, \tau}_{\Delta t},
\end{split}
\end{equation*}
which verifies (1) in Lemma \ref{lem:bootstrap}. 
	
At this point we have verified assertions (1)--(4) of Lemma \ref{lem:bootstrap} 
for our choice of $\vc{H}(t, \tau)$ and $\vc{C}(t,\tau)$. Consequently, Lemma \ref{lem:bootstrap} 
tells us that $\vc{C}(t,\tau)$ is true for all times $(t, \tau) \in \Om_{\Delta t}$. In other 
words,
\begin{equation}\label{eq:final2}
	\|\vvt(t, \tau)\|_{H^{k-2}} \leq \frac{\gamma}{2}, \quad  (t,\tau) \in \Om_{\Delta t},
\end{equation}
which concludes the existence part of Theorem \ref{thm:main}.

Equipped with \eqref{eq:final2}, we can apply Lemma \ref{lem:w} to obtain 
the error estimate
\begin{equation*}
	\|e(t)\|_{H^{k-2}} \leq t \Delta t\, C e^{Ct}\|\vt_0\|^2_{H^k}, \quad t \in [0,T],
\end{equation*}
which is the second part of Theorem \ref{thm:main}.\qed

\section{Strang splitting}\label{sec:strang}
In this section we prove the expected second-order convergence rate of 
Strang splitting applied to \eqref{eq:qg}. 
We  also prove that the method 
 produces regular solutions up to any 
given finite final time $T>0$. The results are 
valid under a condition on the size 
of each time step $\Delta t$.

We will continue to use the notation and definitions introduced in 
the previous section, unless explicitly stating otherwise.
The Strang method we consider in this section reads as
follows: For $\Delta t> 0$ given,  construct a sequence 
$\{\vt^n\}_{n=0}^{\lfloor T/\Delta t\rfloor}$ of approximate solutions 
to \eqref{eq:qg} by the following procedure:
Let $\vt^0 = \vt_0$ and determine sequentially
\begin{equation*}
	\begin{split}
	\vt^{n+1} &= \Phi_B\left(\frac{\Delta t}{2}, 
	\Phi_A(\Delta t, \Phi_B(\frac{\Delta t}{2}, \vt^{n}))\right) \\
	&= \left[\Phi_B(\frac{\Delta t}{2})\circ 
	\Phi_A(\Delta t)\circ \Phi_B(\frac{\Delta t}{2})\right](\vt^n).
	\end{split}
\end{equation*}

As with the Godunov splitting, the convergence analysis will 
require a time-continuous interpolation of the Strang approximations 
$\{\vt^n\}_n$. In contrast to the Godunov case, we will now introduce 
three time variables instead of two. This approach 
is different from the one taken in \cite{Holden:tao}, and indeed appears more 
natural. In \cite{Holden:tao} the authors stick to two time variables and 
interpret Strang splitting in terms of two ``$\Delta t/2$"
Godunov splittings, and in alternating order. 

\begin{figure}[tbp]
  \centering
  \includegraphics[width=0.8\linewidth]{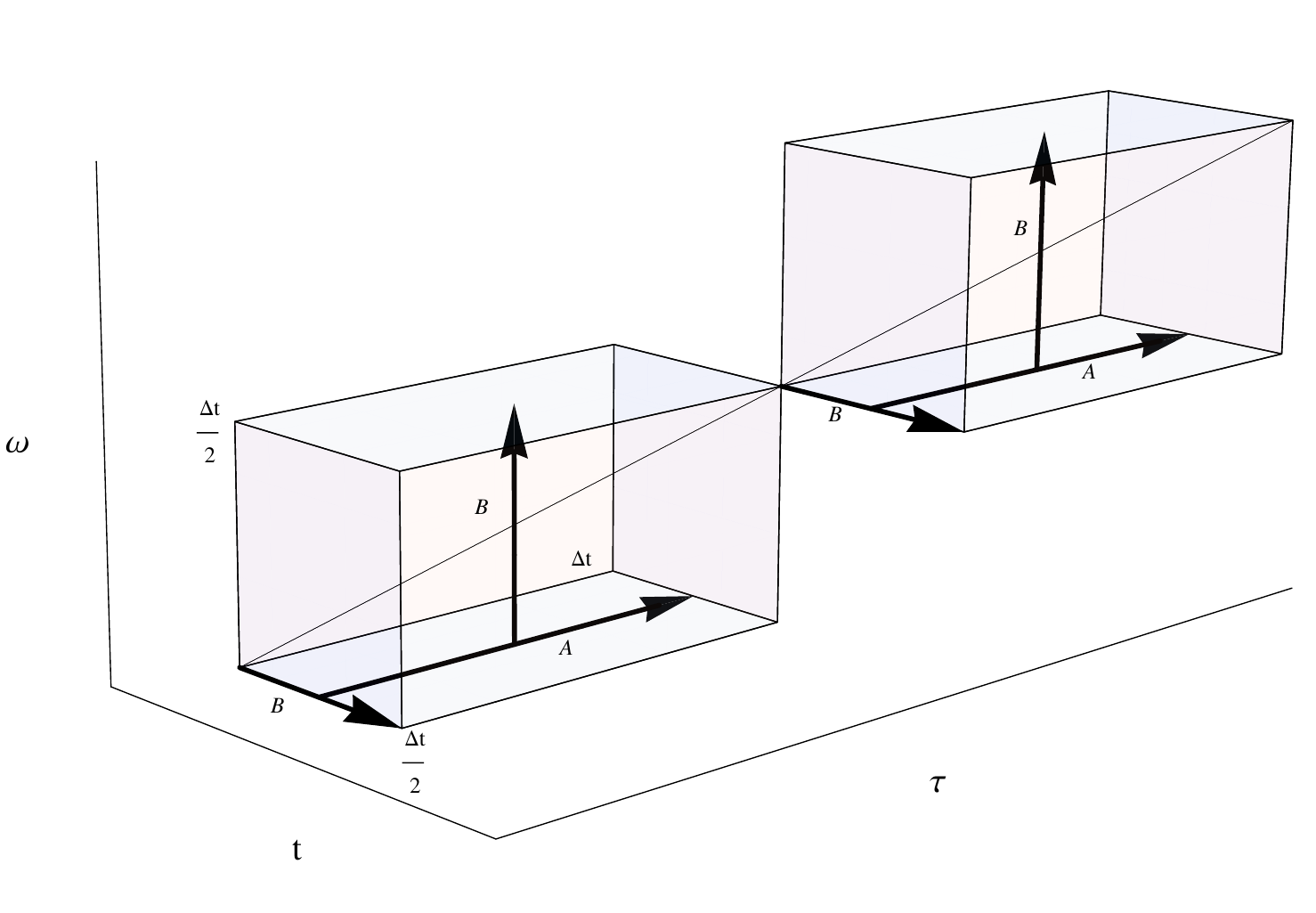}
  \caption{A schematic view of Strang splitting.}
  \label{fig:StrangSplit}
\end{figure}

\begin{definition}[Strang splitting]\label{def:strang}
With $\Delta t > 0$ given, define the domain
\begin{equation*}
	\Om_{\Delta t} =  \bigcup_{n=0}^{\lfloor T/ \Delta t\rfloor-1}
	\left [\frac{t_n}{2}, \frac{t_{n+1}}{2}\right] \times 
	\left [t_n, t_{n+1}\right] \times 
	\left [\frac{t_n}{2}, \frac{t_{n+1}}{2}\right].
\end{equation*}
The time-continuous Strang splitting solution
$\vvt:\Om_{\Delta t} \to \R$ is defined as
\begin{equation}\label{eq:strang-2}
	\begin{split}
		&\vvt(0,0,0) = \vt_0, \\
		& \begin{cases}
			\vvt_t\left(t,t_n,\frac{t_n}{2}\right) = 
			B\left(\vvt\left(t,t_n,\frac{t_n}{2}\right)\right), 
			& t \in \left(\frac{t_n}{2}, \frac{t_{n+1}}{2}\right],\\ 
			\vvt|_{t=\frac{t_n}{2}}=\vt^n;
		 \end{cases}
		 \\ &
		 \begin{cases}
			\vvt_\tau\left(t, \tau,\frac{t_n}{2}\right)
			=A\left(\vvt\left(t, \tau,\frac{t_n}{2}\right)\right), 
			& (t,\tau) \in \left[\frac{t_n}{2}, \frac{t_{n+1}}{2}\right]
			\times \left(t_n,t_{n+1}\right], \\
			\vvt|_{\tau=t_n}=\vvt\left(t,t_n,\frac{t_n}{2}\right), & 
			t\in \left[\frac{t_n}{2}, \frac{t_{n+1}}{2}\right];
		\end{cases}
		\\ &
		\begin{cases} 
			\vvt_\omega(t, \tau, \omega)= B\left(\vvt\left(t,\tau, \omega\right)\right), 
			& (t, \tau, \omega) \in \left[\frac{t_n}{2}, \frac{t_{n+1}}{2}\right]
			\times [t_n,t_{n+1}] \times \left(\frac{t_n}{2}, \frac{t_{n+1}}{2}\right),\\
			\vvt|_{\omega=\frac{t_n}{2}}=\vvt\left(t,\tau,\frac{t_n}{2}\right),
			& (t,\tau) \in \left[\frac{t_n}{2}, \frac{t_{n+1}}{2}\right],
		\end{cases}
	\end{split}
\end{equation}
cf.~Figure \ref{fig:StrangSplit}.
\end{definition}

In each box $\left[\frac{t_n}{2}, \frac{t_{n+1}}{2}\right] 
\times \left[t_n, t_{n+1}\right] \times \left[\frac{t_n}{2}, \frac{t_{n+1}}{2}\right]$, it is 
of particular interest to consider the function $\vvt$ restricted 
to the diagonal, i.e., the function $\vvt\left(\frac{t}{2},t,\frac{t}{2}\right)$ 
for $t \in [t_n, t_{n+1}]$. Observe that at each point on the diagonal, $\vvt$ 
is a Strang splitting solution with a specific time step. 
Specifically, $\vvt\left(\frac{t}{2},t,\frac{t}{2}\right)$ is a 
Strang splitting solution with time step $t - t_n$. 
Consequently,
\begin{equation*}
	\vvt\left(\frac{t_n}{2}, t_n,\frac{t_n}{2}\right) 
	= \vt^n, \quad n=0, \ldots, \lfloor T/\Delta t\rfloor,
\end{equation*}
and hence that $\vvt\left(\frac{t}{2},t,\frac{t}{2}\right)$ can be 
seen as an extension of $\{\vt^n\}_n$ to all of $[0,T]$.

To measure the error between the splitting 
approximation and the exact solution, we will use the function
\begin{equation*}
	e(t) =  \vvt\left(\frac{t}{2},t,\frac{t}{2}\right)-\vt(t),
\end{equation*}
where $\vt$ is the (smooth) solution of \eqref{eq:qg}.

We state two main results in this section. 

\begin{lemma}[Well-defined]\label{thm:strang}
Let $\vvt=\vvt(t,\tau,\omega)$ be the Strang splitting solution of \eqref{eq:qg}
in the sense of  Definition \ref{def:strang}. Suppose 
$\theta_0 \in H^k$ with $6 \leq k \in \mathbb{N}$, $\alpha$, $\beta$ $\in [1,2]$, and that 
$\Delta t>0$ is sufficiently small. Then
\begin{enumerate}
	\item the Strang splitting method  
	is well-defined with $\vvt \in C([0,T]^2;H^k)$;
	\item the error $e(t)=\vvt(\frac{t}{2},t,\frac{t}{2}) - \vt(t)$ satisfies
	$$
	\|e(t)\|_{H^{k-2}} \leq C\, t\, \Delta t.
	$$
	for some constant $C=C\left(\|\theta_0\|_{H^k}\right)$.
\end{enumerate}
\end{lemma}

\begin{theorem}[Convergence]\label{thm:rate}
Under the conditions of the previous lemma,
\begin{equation*}
	\|e(t)\|_{H^{k-3\alpha }} \leq C(\Delta t)^2,
\end{equation*}
for some constant $C=C\left(\|\theta_0\|_{H^k}\right)$.
\end{theorem}
Lemma \ref{thm:strang} and Theorem \ref{thm:rate} will 
be consequences of the results stated and proved in 
the ensuing subsections. 

\subsection{The error evolution equations}
To prove Lemma \ref{thm:strang} we will use 
the same approach as for the Godunov method. However, 
since we are now using three time variables instead of 
two, we must derive a new set of evolution equations governing the error.

By a direct calculation, the error $e$ satisfies the time-evolution
\begin{equation*}
	\begin{split}
		e_t - dC(\vt)[e] 
		&= \frac{\vvt_t}{2} + \vvt_\tau  + \frac{\vvt_\omega}{2} 
		- \vt_t - dA(\vt)[e] - dB(\vt)[e] \\
		&= \frac{\vvt_t}{2} + \vvt_\tau + \frac{1}{2}B(\vvt) - (A+B)(\vt) 
		- dA(\vt)[e] - dB(\vt)[e] \\
		&= \frac{1}{2}\left(\vvt_t - B(\vvt)\right) + \vvt_\tau - A(\vvt) \\
		&\qquad + (A(\vvt) - A(\vt) - dA(\vt)[e]) 
		 + (B(\vvt) - B(\vt) - dB(\vt)[e])  \\
		&= F(t) + \int_0^1 (1-\gamma)d^2 C(\vt 
		+ \gamma e)[e]^2~d\gamma,
	\end{split}
\end{equation*}
where $F(t) = F\left(\frac{t}{2}, t, \frac{t}{2}\right)$ and
\begin{equation*}
	F(t,\tau,\omega) = \frac{1}{2}\left(\vvt_t(t,\tau,\omega) 
	- B(\vvt(t,\tau, \omega))\right)
	+ \vvt_\tau - A(\vvt(t,\tau,\omega)).
\end{equation*}

Since $\vvt_\tau - A(\vvt(t,\tau,\omega)) = 0$ when $\omega =
\frac{t_n}{2}$, $n=0, \ldots,\lfloor T/\Delta t\rfloor-1$, we 
can apply the arguments of \eqref{eq:err3} to obtain 
\begin{equation*}
	F_\tau - dA(\vvt)[F] 
	= \frac{1}{2}[A,B](\vvt), \quad \left(t,\tau,\frac{t_n}{2}\right) \in \Om_{\Delta t},
\end{equation*}
where 
\begin{equation*}
	[A,B](f) = dA(f)[B(f)] - dB(f)[A(f)].
\end{equation*}

We also derive the following equation for the evolution of $F$ in $\omega$:
\begin{equation*}
	\begin{split}
		&F_\omega  - dB(\vvt)[F] \\
		&\qquad= \frac{1}{2}\vvt_{t\omega} - \frac{1}{2}B(\vvt)_\omega 
		- \frac{1}{2}dB(\vvt)[\vvt_t - B(\vvt)] \\
		&\qquad \qquad + \vvt_{\tau \omega} - dA(\vvt)[\vvt_\omega] 
		- dB(\vvt)[\vvt_\tau - A(\vvt)]\\
		&\qquad= \frac{1}{2}B(\vvt)_t - \frac{1}{2}dB(\vvt)[\vvt_\omega] 
		- \frac{1}{2}dB(\vvt)[\vvt_t] +  \frac{1}{2}dB(\vvt)[B(\vvt)] \\
		&\qquad \qquad + dB(\vvt)[\vvt_\tau] - dA(\vvt)[B(\vvt)] - dB(\vvt)[\vvt_\tau] 
		+ dB(\vvt)[A(\vvt)] \\
		&\qquad= \frac{1}{2}\left(dB(\vvt)[\vvt_t] - dB(\vvt)[B(\vvt)]- dB(\vvt)[\vvt_t] 
		+  dB(\vvt)[B(\vvt)]\right) \\
		&\qquad \qquad + dB(\vvt)[A(\vvt)] - dA(\vvt)[B(\vvt)] \\
		&\qquad =dB(\vvt)[A(\vvt)] - dA(\vvt)[B(\vvt)].
	\end{split}
\end{equation*}
Thus, recalling the definition of $[\dott,\dott]$, 
$$
F_\omega  - dB(\vvt)[F] = [B,A](\vvt).
$$
For \eqref{eq:qg} the above evolution equations take the form:
\begin{align}
	F &= \frac{1}{2}\left(\vvt_t 
	+ \Grad \vvt\cdot \Curl \Lambda^{-\beta}\vvt\right)
	+ \vvt_\tau + \Lambda^\alpha \vvt, 
	\quad  (t,\tau, \omega) \in \Om_{\Delta t}, \notag 
	\\ F_\tau + \Lambda^\alpha F 
	&= \frac{1}{2}G^\alpha(\Grad \vvt, \Curl \Lambda^{-\beta}\vvt), 
	\quad \omega = \frac{t_n}{2}, \, n=0,\ldots, \lfloor T/\Delta t\rfloor, \label{eq:fplane}\\
	F_\omega &= - \Grad \vvt\cdot \Curl \Lambda^{-\beta} F 
	-\Grad F\cdot \Curl \Lambda^{-\beta}\vvt \nonumber\\
	& \qquad -G^\alpha(\Grad \vvt, \Curl \Lambda^{-\beta}\vvt), 
	\quad  (t,\tau, \omega) \in \Om_{\Delta t}, \label{eq:fomega}\\
	e_t + \Lambda^\alpha e &= 
	F - \Grad e\cdot \Curl \Lambda^{-\beta}\vt 
	- \Grad \vt\cdot \Curl \Lambda^{-\beta}e \nonumber\\
	&\qquad - \Grad e\cdot \Curl \Lambda^{-\beta}e, 
	\quad  (t,\tau, \omega) \in \Om_{\Delta t}.\label{eq:erroreq}
\end{align}

To prove Lemma \ref{thm:strang} we will adapt 
Lemmas \ref{lem:higher}, \ref{lem:F}, and \ref{lem:w} to the 
Strang splitting approximation. 

Let $\Om_{\Delta t}^{t,\tau,\omega}$  denote the set 
of times prior to $(t, \tau, \omega)$:
\begin{equation}\label{eq:recall}
	\Om_{\Delta t}^{t,\tau, \omega}
	=\Om_{\Delta t}\cap\big(\,[0,t]\times[0,\tau]\times[0,\omega]\, \big).
\end{equation}

We begin by observing that Lemma \ref{lem:higher} can be
easily adapted to yield the following result.

\begin{lemma}\label{lem:higher2}
Let $6 \leq k \in \mathbb{N}$, and assume the 
existence of $(t,\tau,\omega) \in \Om_{\Delta t}$ such that
\begin{equation}
	\|\vvt(s,\sigma, \zeta)\|_{H^{k-2}} \leq \gamma, 
	\qquad  (s,\sigma,\zeta) \in \Om_{\Delta t}^{t,\tau,\omega}.
\end{equation}
Then
\begin{equation}
	\|\vvt(s,\sigma, \zeta)\|_{H^k} \leq e^{C\gamma (s+\zeta)}
	\|\theta_0\|_{H^k},\qquad  (s,\sigma,\zeta) 
	\in \Om_{\Delta t}^{t,\tau,\omega}.
\end{equation}
\end{lemma}

Next we prove a Strang version of Lemma \ref{lem:F}.

\begin{lemma}\label{lem:F2}
Under the conditions of the previous lemma,
\begin{equation*}
	\|F(s,\sigma,\zeta)\|_{H^{k-2}} \leq 
	\Delta t\, C_1 e^{C_2\gamma(s+\zeta)}, 
	\qquad  (s, \sigma, \zeta) \in \Om_{\Delta t}^{t,\tau, \omega}.
\end{equation*}
\end{lemma}

\begin{proof}

First, we observe that in the plane given by $\omega = \frac{t_n}{2}$,
\begin{equation*}
	F(s,\sigma, \frac{t_n}{2}) = \frac{1}{2}\left(\vvt_t - B(\vvt)
        \right), \quad n=0, \ldots, \lfloor T/\Delta t\rfloor-1.
\end{equation*}
Note that while $F \neq 0$ on the line $(s, t_n, \frac{t_n}{2})$,
there holds 
 $\frac{1}{2}\left(\vvt_t - B(\vvt) \right) = 0$ on $(s, t_n, \frac{t_n}{2})$.
Hence, we can repeat the arguments of Lemma \ref{lem:F} to conclude that
\begin{equation}\label{eq:est1}
	\begin{split}
		\|F(s,\sigma, \frac{t_n}{2})\|_{H^{k-2}} 
		&= \|\frac{1}{2}\left(\vvt_t(s,\sigma, \frac{t_n}{2}) - B(\vvt(s,\sigma, \frac{t_n}{2}))
	        \right)\|_{H^{k-2}} \\
		&\leq \Delta t\, C_1
	        e^{C_2\gamma s}, \quad  (s,\sigma, \frac{t_n}{2}) \in \Omega_{\Delta t}^{t, \tau, \omega}.		
	\end{split}
\end{equation}

We fix $n$  such that
$(t,\tau,\omega)\in[\frac{t_n}{2},\frac{t_{n+1}}{2}]\times[t_n,t_{n+1}]\times[\frac{t_n}{2},\frac{t_{n+1}}{2}]$. 
Now, to estimate $\|F\|_{H^{k-2}}$ at an arbitrary point $(s,\sigma, \zeta) \in \Omega_{\Delta t}^{t, \tau, \omega}$, 
we will first integrate  in the $\omega$ direction to the plane given
by $\omega= \frac{t_n}{2}$ and then apply the estimate \eqref{eq:est1}. 
To perform this integration, we apply $\Grad^{s}$ to \eqref{eq:fomega}, 
multiply the result with $\Grad^{s}F$, sum over $s=0,\ldots, k-2$, and 
integrate by parts, to obtain
\begin{equation}\label{eq:est2}
	\begin{split}
		\partial_\omega \frac{1}{2}\|F\|_{H^{k-2}}^2 
		&= -\sum_{s=0}^{k-2}\int_{\R^2} 
		\Grad^{s}F:\Grad^{s}G^\alpha(\Grad \vvt,\Curl \Lambda^{-\beta} \vvt)~dx \\
		&\qquad- \sum_{s=0}^{k-2}\int_{\R^2} 
		\Grad^{s}\left(\Grad \vvt\cdot\Curl \Lambda^{-\beta} F \right):\Grad^{s}F~dx \\
		&\qquad- \sum_{s=0}^{k-2}\int_{\R^2} 
		\Grad^{s}\left(\Grad F\cdot\Curl \Lambda^{-\beta} \vvt \right):\Grad^{s}F~dx \\
		&:= I_1 + I_2 + I_3,
	\end{split}
\end{equation}
and the goal is to bound  the terms $I_1$, $I_2$, and $I_3$. 

By first estimating as in \eqref{eq:F1}, then using
Corollary \ref{cor:com}, and finally applying Lemma \ref{lem:higher2}, we obtain
\begin{align*}
	|I_1| &\leq  \|F\|_{H^{k-2}}
	\|G^\alpha(\Grad \vvt, \Curl \Lambda^{-\beta} \vvt)\|_{H^{k-2}} \\
	& \le C \|F\|_{H^{k-2}} \|\Grad \vvt\|_{H^{k-1}} 
	\|\Grad  \Curl \Lambda^{-\beta} \vvt\|_{H^{k-1}} \\
	& \le C_1\|F\|_{H^{k-2}} e^{C_2 \gamma (s+ \zeta)} 
	\|\theta_0\|_{H^k}^2. 
\end{align*} 

To bound the two other terms $I_2$ and $I_3$, 
we first apply Lemma \ref{lem:lineardiv} to find
\begin{equation*}
	\begin{split}
		|I_2 + I_3| \leq C\|F\|_{H^{k-2}}^2\|\vvt\|_{H^{k-1}}.
	\end{split}
\end{equation*}
An application of Lemma \ref{lem:higher2}  then yields
\begin{equation*}
	\begin{split}
		|I_2 + I_3|
		\leq  C \|F\|_{H^{k-2}}^2 
		e^{C_2 \gamma (s+ \zeta)} \|\theta_0\|_{H^k}. 
	\end{split}
\end{equation*}
Inserting these bounds back into \eqref{eq:est2},
\begin{equation*}
	\partial_\omega \|F\|_{H^{k-2}} 
	\leq C_1 e^{C_2 \gamma (s+ \zeta)} + C_3 \|F\|_{H^{k-2}}.
\end{equation*}
Applying the Gronwall lemma we obtain
\begin{equation*}
	\|F(s,\sigma, \zeta)\|_{H^{k-2}} \leq 
	e^{C_3\gamma \Delta t}\left(\left\|F\left(s,\sigma, \frac{t_n}{2}\right)\right\|_{H^{k-2}} 
	+ \Delta t\, C_1 e^{C_2 \gamma (s+ \zeta)}\right).
\end{equation*}
Finally, we apply \eqref{eq:est1} to conclude 
\begin{equation*}
	\|F(s,\sigma, \zeta)\|_{H^{k-2}}\leq \Delta t\, C_1 e^{C_2 \gamma (s+ \zeta)},
\end{equation*}
which completes the proof.
\end{proof}

\begin{lemma}\label{lem:W2}
Under the conditions of Lemma \ref{lem:higher2},
\begin{equation*}
	\|e(s)\|_{H^{k-2}} \leq s \Delta t\, C(\gamma), \quad s \in [0,t],
\end{equation*}	
where $C(\gamma)$ is a constant depending only on $\gamma$, $T$, and 
$\|\vt_0\|_{H^k}$.
\end{lemma}

\begin{proof}
Since \eqref{eq:erroreq} is identical to \eqref{eq:weq}, we can repeat 
the proof of Lemma \ref{lem:w} step by step (using 
Lemmas \ref{lem:higher2} and \ref{lem:F2} 
instead of Lemmas \ref{lem:higher} and \ref{lem:F}) to conclude.	
\end{proof}

\subsection{Proof of Lemma \ref{thm:strang}}
Recall that  $\Om_{\Delta t}^{t,\tau, \omega}$ denotes the set of all
times prior to $(t, \tau, \omega)$, cf.~\eqref{eq:recall}.
To prove the existence part of Lemma \ref{thm:strang} we will utilize
the same strategy as with did for the Godunov method in the previous section. That is, 
to apply the bootstrap lemma (Lemma \ref{lem:bootstrap}) with
 $\vc{H}(t, \tau, \omega)$ as the statement
\begin{equation}\label{eq:Hstrang}
	\|\vvt(s,\sigma, \zeta)\|_{H^{k-2}} \leq \gamma, 
	\quad  (s,\sigma, \zeta) \in \Om^{t,\tau, \omega}_{\Delta t},
\end{equation}
and $\vc{C}(t, \tau, \omega)$  the statement
\begin{equation*}
	\|\vvt(s,\sigma, \zeta)\|_{H^{k-2}} \leq \frac{\gamma}{2}, 
	\quad  (s,\sigma, \zeta) \in \Om^{t,\tau, \omega}_{\Delta t},
\end{equation*}
where $\gamma$ is some value to be specified below. 
Note that there is no problem with extending Lemma \ref{lem:bootstrap} so that it 
allows three time variables instead of two (i.e., hypothesis 
and statement depending on three variables).

To conclude the existence part of Lemma \ref{thm:strang} we must 
verify assertions (1)--(4) in Lemma \ref{lem:bootstrap}.
Assertions  (2) and (3) clearly hold. Regarding assertion (4), we 
make the assumption
\begin{equation}\label{eq:asmpt2}
	\|\vvt(0,0,0)\|_{H^{k-2}} = \|\vt_0\|_{H^{k-2}} \leq \gamma.
\end{equation}

It remains to verify assertion (1). For this purpose, let us 
assume that $\vc{H}(t, \tau, \omega)$ is true 
for some $(t, \tau, \omega) \in \Om_{\Delta t}$.
In order to show assertion (1) we need to prove 
that $\vc{C}(t,\tau,\omega)$ is  true. We first note 
that $\vc{H}(t,\tau,\omega)$ being true
renders Lemma \ref{lem:W2} applicable, allowing us to conclude
\begin{equation}\label{eq:vf12}
	\|e(s)\|_{H^{k-2}} \leq s \Delta t\, C_1(\gamma), \quad s \in [0,t].
\end{equation}

Now, to estimate $\|\vvt(s,\sigma, \zeta)\|_{H^{k-2}}$ at 
any time $(s,\sigma, \zeta) \in \Om_{\Delta t}^{t, \tau, \omega}$ 
the idea is to consider the time evolution of $\|\vvt(s,\sigma, \zeta)\|_{H^{k-2}}$ in 
the $\omega$ and $\tau$ directions. Let $n$ be such that  $(t,\tau,\omega)
\in[\frac{t_n}{2},\frac{t_{n+1}}{2}]\times[t_n,t_{n+1}]\times[\frac{t_n}{2},\frac{t_{n+1}}{2}]$.
First we integrate from $(s,\sigma, \zeta)$ to $(s, \sigma, \frac{t_n}{2})$. 
Then, we integrate the result in the plane 
$\omega= \frac{t_n}{2}$ from $(s,\sigma, \frac{t_n}{2})$ to $(s, 2s,  \frac{t_n}{2})$. Finally, 
we integrate in the $\omega$ direction from $(s, 2s,  \frac{t_n}{2})$ to the diagonal $(s,2s,s)$.
On the diagonal, we will then utilize \eqref{eq:vf12} to 
conclude the necessary bound. 
	
By definition, 
\begin{equation}\label{eq:lm1}
	\begin{split}
		\left|\partial_\omega \frac{1}{2}\|\vvt (s,\sigma, \zeta)\|_{H^{k-2}}^2\right|
		&= \left|\sum_{s=0}^{k-2}\int_{\R^2} \Grad^{s} \vvt_\omega: \Grad^{s} \vvt~dx\right| \\
		&= \left|-\sum_{s=0}^{k-2}\int_{\R^2} 
		\Grad^{s} \left(\Grad \vvt\cdot \Curl \Lambda^{-\beta}\vvt \right): \Grad^{s} \vvt~dx \right|\\
		&\leq C\|\vvt\|_{H^{k-4}}\|\vvt\|_{H^{k-2}}^2 \leq \gamma C\|\vvt(s,\sigma, \zeta)\|_{H^{k-2}}^2,
	\end{split}
\end{equation}
where the estimate follows from Lemma \ref{lem:nonlineardiv}, and where we
have used \eqref{eq:Hstrang}. 
The fundamental theorem of calculus and \eqref{eq:lm1} 
provides us with the estimate
\begin{equation}\label{eq:lm2}
	\begin{split}
		\|\vvt (s,\sigma, \zeta)\|_{H^{k-2}} 
		&\leq \|\vvt (s,\sigma, \zeta_2)\|_{H^{k-2}} 
		 + C\Delta t \sup_{\zeta' \in [\frac{t_n}{2},\frac{t_{n+1}}{2}]}
		 \|\vvt(s,\sigma,\zeta')\|_{H^{k-2}} \\
		&\leq \|\vvt (s,\sigma, \zeta_2)\|_{H^{k-2}} 
		 + C\Delta t\, \|\vt_0\|_{H^{k-2}}, 		
	\end{split}
\end{equation}
for any $\frac{t_n}{2} \leq \zeta_2 \leq  \frac{t_{n+1}}{2}$.

Similarly, 
\begin{equation*}
	\begin{split}
		\left|\partial_\tau \frac{1}{2}
		\left\|\vvt\left(s,\sigma, \frac{t_n}{2}\right)\right\|_{H^{k-2}}^2 \right|
		&= \left|\sum_{s=0}^{k-2}
		\int_{\R^2} \Grad^{s} \vvt_\tau: \Grad^{s} \vvt~dx\right| \\
		&= \left|-\sum_{s=0}^{k-2} \int_{\R^2}
		\Grad^{s}\Lambda^{\alpha}\vvt:\Grad^{s}\vvt~dx\right| \\
		&\leq \left\|\vvt\left(s,\sigma,  \frac{t_n}{2}\right)\right\|_{H^{k-2}}
		\left\|\vvt\left(s, \sigma,  \frac{t_n}{2}\right)\right\|_{H^{k}},
	\end{split}
\end{equation*}
where we have used that $\alpha \leq 2$. Thus
\begin{equation*}
	\left|\partial_\tau \frac{1}{2}
	\left\|\vvt\left(s,\sigma, \frac{t_n}{2}\right)\right\|_{H^{k-2}} \right|
	\le \left\|\vvt\left(s, \sigma,  \frac{t_n}{2}\right)\right\|_{H^{k}}.
\end{equation*}
From the fundamental theorem of calculus, using the previous inequality 
 and  Lemma \ref{lem:higher2}, we conclude 
\begin{equation}\label{eq:lm3}
	\begin{split}
		\left\|\vvt\left(s, \sigma, \frac{t_n}{2}\right)\right\|_{H^{k-2}} 
		&\leq \left\|\vvt\left(s, 2s, \frac{t_n}{2}\right)\right\|_{H^{k-2}}
		\\ & \qquad 
		+ \Delta t \sup_{\sigma' \in [\min\{2s,\sigma\},\max\{2s,\sigma\}]}
		\left\|\vvt\left(s,\sigma', \frac{t_n}{2}\right)\right\|_{H^k} \\
		& \leq \left\|\vvt\left(s, 2s, \frac{t_n}{2}\right)\right\|_{H^{k-2}}
		+ \Delta t\, C\|\theta_0\|_{H^{k}}.
	\end{split}
\end{equation}
Thus, by combining \eqref{eq:lm2} and \eqref{eq:lm3},
\begin{equation}\label{eq:prev}
	\begin{split}
		\|\vvt (s,\sigma, \zeta)\|_{H^{k-2}} 
		&\leq \left\|\vvt\left(s,\sigma, \frac{t_n}{2}\right)\right\|_{H^{k-2}} 
		 + C\Delta t\, \|\vt_0\|_{H^{k-2}}\\
		&\leq \left\|\vvt\left(s, 2s,\frac{t_n}{2}\right)\right\|_{H^{k-2}}
		+ \Delta t\, C\,\|\theta_0\|_{H^{k}}  \\
		&\leq \left\|\vvt\left(s,2s, s\right)\right\|_{H^{k-2}}
		+ \Delta t\, C\,\|\theta_0\|_{H^{k}} \\
		&\leq \|e(2s)\|_{H^{k-2}} + \|\theta(2s)\|_{H^{k-2}} 
		+\Delta t\, C\,\|\theta_0\|_{H^{k}}\\
		& \leq  C_1(\gamma)\Delta t + C_2,
	\end{split}
\end{equation}
where the last inequality is \eqref{eq:vf12}.
Now, we fix $\gamma$ and $\Delta t$ according to
\begin{equation*}
	\gamma = 4C_2, \qquad \Delta t \leq \frac{C_2}{C_1(\gamma)},
\end{equation*}
and note that this is not conflict with \eqref{eq:asmpt2}. Then, \eqref{eq:prev} gives
\begin{equation*}
	\begin{split}
		\|\vvt (s,\sigma, \zeta)\|_{H^{k-2}}  \leq 2C_2 = \frac{\gamma}{2}.
	\end{split}
\end{equation*}
Since $(s,\sigma, \zeta) \in \Om_{\Delta t}^{t,\tau,\omega}$ was arbitrarily chosen, this verifies $(1)$ in Lemma \ref{lem:bootstrap}. 
We have now verified assertions (1)--(4) of Lemma \ref{lem:bootstrap} 
for $\vc{H}(t, \tau, \omega)$ and $\vc{C}(t,\tau, \omega)$. Consequently, Lemma \ref{lem:bootstrap} 
tells us that $\vc{C}(t,\tau, \omega)$ is true for all times $(t, \tau, \omega) \in \Om_{\Delta t}$. In other 
words,
\begin{equation}\label{eq:final22}
	\|\vvt(t, \tau, \omega)\|_{H^{k-2}} \leq \frac{\gamma}{2}, \quad  (t,\tau, \omega) \in \Om_{\Delta t}.
\end{equation}
Thus, we can conclude the existence part of Lemma \ref{thm:strang}.
Since \eqref{eq:final22} holds, we can apply Lemma \ref{lem:w} to conclude
the error estimate
\begin{equation*}
	\|e(t)\|_{H^{k-2}} \leq t \Delta t\, C, \quad t \in [0,T],
\end{equation*}
which is the second part of Lemma \ref{thm:strang}.

\qed

%
%
%
%

\subsection{Temporal regularity and auxiliary estimates}
So far we have focused on the spatial regularity 
of the splitting solution. To  prove 
Theorem \ref{thm:rate}, we will need regularity 
in time of $\vvt$. Since we have introduced 
three time variables, it is convenient 
to define a time gradient. We will use the notation
\begin{equation*}
	\nabla_{ t}f = 
	\begin{pmatrix}
	 f_t \\	 f_\tau \\ f_\omega
	\end{pmatrix}
	, \qquad \nabla_{ t}^2 f =
	\begin{pmatrix}
		f_{tt} & f_{t \tau} & f_{t \omega} \\
		f_{\tau t} & f_{\tau \tau} & f_{\tau \omega} \\
		f_{\omega t} & f_{\omega \tau} & f_{\omega \omega} 
	\end{pmatrix},
\end{equation*}
with the obvious extension to higher order (i.e., $\nabla_{ t}^k f$
for $k$th derivative).

\begin{lemma}\label{lem:tempv}
Let $\alpha \in [1,2]$ and
 $\vvt$ be the Strang splitting solution in 
the sense of Definition \ref{def:strang}.
Then, for any $k$  
and $l$ such that $\vt_0 \in H^{k+\alpha l}$,
\begin{equation*}
	\|\nabla^l_{ t} \vvt\|_{H^k} \leq C\left(\|\vt_0\|_{H^{k+\alpha
            l}}\right)\,,
\quad (t,\tau,\omega)\in\Om_{\Delta t}.
\end{equation*}
\end{lemma}
\begin{proof}
We argue by induction on $l$. For $l=0$, the result is merely an application of Lemma \ref{lem:higher2}.
Now, let us assume that the bound holds for $l=0, \ldots, q$. To close the induction argument, 
it remains to show the bound for $l=q+1$.

Let $(t',\tau',\omega') \in \Omega_{\Delta t}$ be arbitrary and fix $n$ such 
that 
$$
(t',\tau',\omega') \in \left[\frac{t_n}{2},\frac{t_{n+1}}{2}\right]
\times \left[t_n,t_{n+1}\right]\times\left[\frac{t_n}{2},\frac{t_{n+1}}{2}\right]. 
$$
An arbitrary component of $\nabla_{t}^{q+1} \vvt$ can be written in the form
\begin{equation*}
	\Theta^{q+1}_{i,j,\ell} = \frac{\partial^{q+1}}{\partial t^{i}\partial\tau^{j}
	\partial\omega^{\ell}}\vvt, \quad 0 \leq i,j,\ell \leq q+1, 
	\quad i+j+\ell = q+1, \quad i,j,\ell \in \mathbb{N}.
\end{equation*}
To estimate the arbitrary component at the point $(t',\tau',\omega')$, we first apply
the fundamental theorem of calculus to write
\begin{equation*}
	\begin{split}
		\frac{1}{2}\|\Theta_{i,j,\ell}^{q+1}(t',\tau',\omega')\|_{H^k}^2 
		&
		=\frac{1}{2}\left\|\Theta_{i,j,\ell}^{q+1}\left(t',\tau',\frac{t_n}{2}\right)\right\|_{H^k}^2 \\
		 &\quad + \sum_{s=0}^k\int_{\frac{t_n}{2}}^{\omega'}\!\!
		\int_{\R^2}\Grad^s\partial_\omega\Theta_{i,j,\ell}^{q+1}(t',\tau',s):
		\Grad^s\Theta_{i,j,\ell}^{q+1}(t',\tau',s)~dxds.
	\end{split}
\end{equation*}
By applying the definition of $\Theta_{i,j,\ell}^{q+1}$,
the definition $\vvt_\omega$ (cf.~\eqref{eq:strang-2}), and the Leibniz rule,
we deduce
\begin{equation*}
	\begin{split}
		&\int_{\R^2}\Grad^s\partial_\omega\Theta_{i,j,\ell}^{q+1}(t',\tau',s):
		\Grad^s\Theta_{i,j,\ell}^{q+1}(t',\tau',s)~dxds \\
		&\quad = 
		\sum_{s=0}^k\int_{\frac{t_n}{2}}^{\omega'}\!\!
				\int_{\R^2}\Grad^s\frac{\partial^{q+1}}{\partial t^{i}\partial\tau^{j}
				\partial\omega^{\ell}}\vvt_\omega(t', \tau', s):
				\Grad^s\Theta_{i,j,\ell}^{q+1}(t',\tau',s)~dxds \\
		&\quad = 
		\sum_{s=0}^k \sum_{r=0}^i \sum_{n=0}^j\sum_{m=0}^{\ell}
		\left(i \atop r\right)\left(j \atop n\right)\left(\ell \atop m\right)\\
		&\quad \qquad  \qquad \times \int_{\frac{t_n}{2}}^{\omega'}
		\int_{\R^2}\Grad^s\left(\Curl \Lambda^{-\beta}\Theta_{r,n,m}^{m+n+r}(t',\tau',s)
		\cdot \Grad \Theta_{i-r,j-n,\ell-m}^{l-r-n-m}(t',\tau',s)\right) \\
		&\qquad \qquad \qquad \qquad \qquad \qquad \qquad
		: \Grad^s\Theta_{i,j,\ell}^{q+1}(t',\tau',s)~dxds.
	\end{split}
\end{equation*}
Inserting this expression into the previous equality yields
\begin{align}\label{eq:fundthm}
		&\frac{1}{2}\|\Theta_{i,j,\ell}^{q+1}(t',\tau',\omega')\|_{H^k}^2\nonumber \\
		&\quad = \frac{1}{2}
		\left\|\Theta_{i,j,\ell}^{q+1}\left(t',\tau',\frac{t_n}{2}\right)\right\|_{H^k}^2
		\\ &\quad \qquad + 
		\sum_{s=0}^k \sum_{r=0}^i \sum_{n=0}^j\sum_{m=0}^{\ell}
		\left(i \atop r\right)\left(j \atop n\right)\left(\ell \atop m\right) \nonumber\\
		&\quad \qquad  \qquad \times \int_{\frac{t_n}{2}}^{\omega'}
		\int_{\R^2}\Grad^s\left(\Curl \Lambda^{-\beta}\Theta_{r,n,m}^{m+n+r}(t',\tau',s)
		\cdot \Grad \Theta_{i-r,j-n,\ell-m}^{l-r-n-m}(t',\tau',s)\right) \nonumber\\
		&\qquad \qquad \qquad \qquad \qquad \qquad \qquad
		: \Grad^s\Theta_{i,j,\ell}^{q+1}(t',\tau',s)~dxds. \nonumber
\end{align}

Let us consider three separate cases 
of $m+n+r$ in the quadruple sum above.\\
(i) If $m+n+r = q+1$, the corresponding term in the above reads
\begin{align}
	&\sum_{s=0}^k\int_{\frac{t_n}{2}}^{\omega'}
	\int_{\R^2}\Grad^s\left(\Curl \Lambda^{-\beta}\Theta_{i,j,\ell}^{q+1}(t',\tau',s)
	\cdot \Grad \vvt (t',\tau',s)\right): \Grad^s\Theta_{i,j,\ell}^{q+1}(t',\tau',s)~dxds \nonumber\\
	&\qquad \leq C\int_{\frac{t_n}{2}}^{\omega'}\|\vvt(t',\tau',s)\|_{H^{k+1}}
	\|\Theta_{i,j,\ell}^{q+1}(t',\tau',s)\|_{H^k}^2~ds \nonumber \\
	&\qquad \leq C\|\vt_0\|_{H^{k+1}}\int_{\frac{t_n}{2}}^{\omega'}
	\|\Theta_{i,j,\ell}^{q+1}(t',\tau',s)\|_{H^k}^2~ds.
	\label{eq:indu1}
\end{align}
Here, the first inequality is an application of Lemma \ref{lem:lineardiv} and the last inequality is 
an application of Lemma \ref{lem:higher2}. \\
(ii) If $m+n+r=0$, we can also apply Lemma \ref{lem:lineardiv} to conclude 
\begin{align*}
	&\sum_{s=0}^k\int_{\frac{t_n}{2}}^{\omega'}
	\int_{\R^2}\Grad^s\left(\Curl \Lambda^{-\beta}\vvt (t',\tau',s)
	\cdot \Grad \Theta_{i,j,\ell}^{q+1}(t',\tau',s)\right):
	\Grad^s\Theta_{i,j,\ell}^{q+1}(t',\tau',s)~dxds \\
	&\qquad \leq C\int_{\frac{t_n}{2}}^{\omega'}
	\|\vvt(t',\tau',s)\|_{H^{k+1}}\|\Theta_{i,j,\ell}^{q+1}(t',\tau',s)\|_{H^k}^2~ds  \\
	&\qquad \leq C\|\vt_0\|_{H^{k+1}}
	\int_{\frac{t_n}{2}}^{\omega'}\|\Theta_{i,j,\ell}^{q+1}(t',\tau',s)\|_{H^k}^2~ds.
\end{align*}
(iii) For the remaining cases ($1 \leq m+n+r \leq q$), we first 
apply the Leibniz rule to the spatial derivatives, then 
we overestimate the terms using the H\"older inequality. In this way, 
we obtain 
\begin{align}
	&\sum_{s=0}^k\int_{\frac{t_n}{2}}^{\omega'}
	\int_{\R^2}\Grad^s\left(\Curl \Lambda^{-\beta}\Theta_{r,n,m}^{m+n+r}(t',\tau',s)
	\cdot \Grad \Theta_{i-r,j-n,\ell-m}^{q+1-r-n-m}(t',\tau',s)\right)\nonumber \\
	&\qquad \qquad \qquad\qquad\qquad 
	: \Grad^s\Theta_{i,j,\ell}^{q+1}(t',\tau',s)~dxds \nonumber \\
	&\leq C \int_{\frac{t_n}{2}}^{\omega'}\|\Theta_{i,j,\ell}^{q+1}(t',\tau',s)\|_{H^k}
	\|\Theta_{r,n,m}^{m+n+r}(t',\tau',s)\|_{H^{k+1}}
	\\ & \qquad \qquad \qquad\qquad\qquad
	\times \|\Theta_{i-r,j-n,\ell-m}^{{q}-r-n-m}(t',\tau',s)\|_{H^{k+1}}~ds \nonumber \\
	&\leq C \int_{\frac{t_n}{2}}^{\omega'}
	\|\Theta_{i,j,\ell}^{q+1}(t',\tau',s)\|_{H^k}\|\vt_0\|_{H^{k+\alpha(q+1)}}^2~ds,
	\label{eq:indu3}
\end{align}
where we have used the induction hypothesis and $\alpha \geq 1$ to derive the  
last inequality.

By applying \eqref{eq:indu1}--\eqref{eq:indu3} to \eqref{eq:fundthm}, we gather
\begin{equation*}
	\begin{split}
		&\frac{1}{2}\|\Theta_{i,j,\ell}^{q+1}(t',\tau',\omega')\|_{H^k}^2 \\
		&\qquad \leq \frac{1}{2}
		\left\|\Theta_{i,j,\ell}^{q+1}\left(t',\tau',\frac{t_n}{2}\right)\right\|_{H^k}^2 
		+ C\|\vt_0\|_{H^{k+(1+q)\alpha}}^2\int_{\frac{t_n}{2}}^{\omega'}
		\|\Theta_{i,j,\ell}^{q+1}(t',\tau',s)\|_{H^k}~ds  \\
		&\qquad \qquad \qquad 
		+ C\|\vt_0\|_{H^{k+1}}\int_{\frac{t_n}{2}}^{\omega'}
		\|\Theta_{i,j,\ell}^{q+1}(t',\tau',s)\|_{H^k}^2~ds \\
		& \qquad \leq \frac{1}{2}\left\|\Theta_{i,j,\ell}^{q+1}\left(t',\tau',\frac{t_n}{2}\right)
		\right\|_{H^k}^2 + C\Delta t + C\int_{\frac{t_n}{2}}^{\omega'}
		\|\Theta_{i,j,\ell}^{q+1}(t',\tau',s)\|_{H^k}^2~ds,
	\end{split}
\end{equation*}
where the constant $C$ depends on $\|\vt_0\|_{H^{k+(1+q)\alpha}}$. 
By applying Gronwall's inequality to the previous inequality, we conclude
\begin{equation}\label{eq:firststep}
	\begin{split}
		\|\Theta_{i,j,\ell}^{q+1}(t',\tau',\omega')\|_{H^k}^2
		\leq \left(\left\|\Theta_{i,j,\ell}^{q+1}\left(t',\tau',\frac{t_n}{2}\right)\right\|_{H^k}^2 
		+ C\Delta t\right) e^{C\Delta t}.
	\end{split}
\end{equation}

We  have now derived a bound on $\Theta_{i,j,\ell}^{q+1}(t',\tau',\omega')$ in terms 
of  $\Theta_{i,j,\ell}^{q+1}(t',\tau',\frac{t_n}{2})$. Next, we derive a bound 
on  $\Theta_{i,j,\ell}^{q+1}(t',\tau',\frac{t_n}{2})$ in terms 
of $\Theta_{i,j,\ell}^{q+1}(t',t_n,\frac{t_n}{2})$. For this purpose, we once 
more apply the fundamental theorem of calculus to obtain 
\begin{equation*}
	\begin{split}
		&\frac{1}{2}\left\|\Theta_{i,j,\ell}^{q+1}\left(t',\tau',\frac{t_n}{2}\right)\right\|_{H^k}^2 \\
		&\quad =\frac{1}{2}\left\|\Theta_{i,j,\ell}^{q+1}\left(t',t_n,\frac{t_n}{2}\right)\right\|_{H^k}^2\\
		&\qquad+\sum_{s=0}^k
		\int_{t_n}^{\tau'}\int_{\R^2}\Grad^s\partial_\tau
		\Theta_{i,j,\ell}^{q+1}\left(t',\tilde s,\frac{t_n}{2}\right):
		\Grad^s\Theta_{i,j,\ell}^{q+1}\left(t',\tilde s,\frac{t_n}{2}\right)~dxd\tilde s \\
		&\quad = \frac{1}{2}\|\Theta_{i,j,\ell}^{q+1}\left(t',t_n,\frac{t_n}{2}\right)\|_{H^k}^2\\
		&\qquad-\sum_{s=0}^k
		\int_{t_n}^{\tau'}\int_{\R^2}\Grad^s \Lambda^\alpha
		\Theta_{i,j,\ell}^{q+1}(t',\tilde s,\frac{t_n}{2}):
		\Grad^s\Theta_{i,j,\ell}^{q+1}\left(t',\tilde s,\frac{t_n}{2}\right)~dxd\tilde s \\
		&\quad = \frac{1}{2}\left\|\Theta_{i,j,\ell}^{q+1}\left(t',t_n,\frac{t_n}{2}\right)\right \|_{H^k}^2
		- \int_{t_n}^{\tau'}
		\left \|\Theta_{i,j,\ell}^{q+1}\left(t',\tilde s,\frac{t_n}{2}\right)
		\right\|_{H^{k+\alpha/2}}~d\tilde s,
	\end{split}
\end{equation*}
where we have also used that $\vvt_\tau = A(\vvt)$ for $\omega= \frac{t_n}{2}$.
It follows that
\begin{equation}\label{eq:secondstep}
	\left\|\Theta_{i,j,\ell}^{q+1}\left(t',\tau',\frac{t_n}{2}\right)\right\|_{H^k}^2 
	\leq \left\|\Theta_{i,j,\ell}^{q+1}\left(t',t_n,\frac{t_n}{2}\right)\right\|_{H^k}^2.
\end{equation}

Finally, we perform our last application 
of the fundamental theorem to obtain
\begin{equation*}\label{eq:fundthm2}
	\begin{split}
		&\frac{1}{2}\left\|\Theta_{i,j,\ell}^{q+1}\left(t',t_n,\frac{t_n}{2}\right)\right\|_{H^k}^2 
		-\frac{1}{2}\left\|\Theta_{i,j,\ell}^{q+1}\left(\frac{t_n}{2},t_n,\frac{t_n}{2}\right)\right\|_{H^k}^2
		\\ &\qquad \qquad
		=\sum_{s=0}^k \int_{\frac{t_n}{2}}^{\omega'}
		\int_{\R^2}\Grad^s\partial_t
		\Theta_{i,j,\ell}^{q+1}\left(\tilde s,t_n,\frac{t_n}{2}\right):
		\Grad^s\Theta_{i,j,\ell}^{q+1}\left(\tilde s,t_n,\frac{t_n}{2}\right)~dxd\tilde s.
	\end{split}
\end{equation*}

By applying the same calculations to the previous 
equations as those used to derive \eqref{eq:firststep}, we arrive at
\begin{equation}\label{eq:laststep}
	\left\|\Theta_{i,j,\ell}^{q+1}\left(t',t_n,\frac{t_n}{2}\right)\right\|_{H^k}^2 
	\leq \left(\left\|\Theta_{i,j,\ell}^{q+1}\left(\frac{t_n}{2},t_n,\frac{t_n}{2}\right)\right\|_{H^k}^2 
	+ C \Delta t\right)e^{C\Delta t}.
\end{equation}

Combining \eqref{eq:firststep}, \eqref{eq:secondstep}, and \eqref{eq:laststep} gives
\begin{equation*}
	\|\Theta_{i,j,\ell}^{q+1}(t',\tau',\omega')\|_{H^k} 
	\leq
        \left\|\Theta_{i,j,\ell}^{q+1}\left(\frac{t_n}{2},t_n,\frac{t_n}{2}\right)
        \right\|_{H^k}e^{2C\Delta t} 
        + C\Delta t\left(e^{2C\Delta t} + e^{C\Delta t} \right),
\end{equation*}
which immediately leads to
\begin{equation}\label{eq:almostdone}
	\begin{split}
		\|\Theta_{i,j,\ell}^{q+1}(t',\tau',\omega')\|_{H^k} 
		&\leq \|\Theta_{i,j,\ell}^{q+1}(0,0,0)\|_{H^k}e^{nC\Delta t}
		+ C\Delta t\sum_{\vr= 1}^Ne^{\vr C\Delta t}		\\
		&\leq \|\Theta_{i,j,\ell}^{q+1}(0,0,0)\|_{H^k}e^{Ct'} + Ct'e^{Ct'}.
	\end{split}
\end{equation}

Finally, let us estimate $\|\Theta_{i,j,\ell}^{q+1}(0,0,0)\|_{H^k}$. By definition, we have that
\begin{equation*}
	\Theta_{i,j,\ell}^{q+1}(0,0,0) =
        \frac{\partial^{q+1}}{\partial t^{i}\partial\tau^{j}\partial\omega^{\ell}}\vvt(0,0,0).
\end{equation*}
At this point, we can apply each of the time-derivatives to $\vvt(0,0,0)$ and use 
the definition \eqref{eq:strang-2} to translate time-derivatives into spatial derivatives. 
Since $\alpha \geq 1$,  it is clear that the case $j=q+1$ contains the 
highest number of spatial derivatives. In this case,
\begin{equation*}
	\partial_\tau^{q+1} \vvt = -\Lambda^\alpha\partial_\tau^{q} \vvt 
	= (-1)^{q+1} \Lambda^{(q+1)\alpha} \vvt.
\end{equation*}
Thus, for any valid combination of $i$,$j$, and $\ell$, we conclude
\begin{equation*}
	\|\Theta_{i,j,\ell}^{q+1}(0,0,0)\|_{H^k} \leq C\|\vt_0\|_{H^{k+\alpha (q+1)}}.
\end{equation*}
From this, and  \eqref{eq:almostdone}, we conclude that the lemma 
holds also for $l=q+1$. 
\end{proof}

The following lemma is almost a corollary of the previous lemma.

\begin{lemma}\label{lem:doubleF}
	Let $\vvt$ be the Strang splitting solution in 
	the sense of Definition \ref{def:strang} and \eqref{eq:strang-2}.
	Then, for any $k$ such that $\|\vt_0\|_{H^{k+3 \alpha}} \leq
        C$, we have 
\begin{equation*}
	\|\nabla_{ t}^2 F\|_{H^k} \leq C\left(\|\vt_0\|_{H^{k+3\alpha}}\right), 
	\quad (t,\tau, \omega) \in \Om_{\Delta t}.
\end{equation*}
\end{lemma}
\begin{proof}
Let $i$ and $j$ denote any one of $t$, $\tau$, or $\omega$. An arbitrary component 
of $\nabla^2_t F$ can then be written $F_{ij}:=\partial_i \partial_j F$. By definition, 
we have that 
\begin{equation*}
	\begin{split}
	F_{ij} &= \partial_i \partial_j \left[	\frac{1}{2}\left(\vvt_t + \Grad \vvt\cdot \Curl
                \Lambda^{-\beta}\vvt\right)
			 + \vvt_\tau + \Lambda^\alpha \vvt\right]\\
	&= \frac{1}{2}\partial_t \vvt_{ij} + \partial_\tau \vvt_{ij} + \frac{1}{2}\Curl \Lambda^{-\beta}\Grad \vvt_{ij} \cdot \Grad \vvt
	 + \frac{1}{2}\Curl \Lambda^{-\beta}\Grad \vvt \cdot \Grad \vvt_{ij} \\
	& \quad + \frac{1}{2}\Curl \Lambda^{-\beta}\Grad \vvt_{j} \cdot \Grad \vvt_i + \frac{1}{2}\Curl \Lambda^{-\beta}\Grad \vvt_i \cdot \Grad \vvt_{j} + \Lambda^\alpha \vvt_{ij}.
	\end{split}
\end{equation*}	
By applying the H\" older inequality together with the previous lemma, we estimate
\begin{equation*}
	\begin{split}
		\|F_{ij}\|_{H^k} 
		& \leq \frac{3}{2}\|\nabla_{ t}^3 \vvt\|_{H^k}
		+ \|\nabla_{t}^2 \vvt\|_{H^{k+\alpha}} \\
		&\qquad +\frac{1}{2}\|\Grad \vvt\cdot \Curl \Lambda^{-\beta}\vvt_{ij}\|_{H^k} 
		+ \frac{1}{2}\|\Grad  \vvt_{ij}\cdot \Curl \Lambda^{-\beta}\vvt\|_{H^k}\\
		& \qquad  + \frac{1}{2}\|\Grad \vvt_i\cdot \Curl \Lambda^{-\beta}\vvt_{j}\|_{H^k}
		+\frac{1}{2}\|\Grad \vvt_j\cdot \Curl \Lambda^{-\beta}\vvt_{i}\|_{H^k}\\
		& \leq C + \|\Grad^{k+1} \vvt\|_{L^\infty}\|\nabla_{ t}^2\vvt\|_{H^k}
		+ \|\nabla_{ t}^2 \vvt\|_{H^k}\|\Grad^k\vvt\|_{L^\infty} \\
		& \qquad + 2\|\nabla_t \vvt\|_{H^k}\|\Lambda^{-\beta}\nabla_t \vvt\|_{H^{k+3}} \\
		& \leq C\left(1+ \|\vvt\|_{H^{k+3}} \right) + C\|\nabla_t \vvt\|_{H^k}\|\nabla_t \vvt\|_{H^{k+2}} \leq C,
	\end{split}
\end{equation*}
where we have grossly overestimated most of the terms.
We have also applied Lemma \ref{lem:tempv}  with $l=1$ for the $H^{k+2}$ norm 
and with $l=2$ for the $H^k$ norm. The constant $C$  
depends on $\|\theta_0\|_{H^{k+3 \alpha}}$.
\end{proof}

Theorem \ref{thm:rate} is
 a  consequence of the following (remarkable) fact:
\begin{lemma}\label{lem:magic}
There holds
\begin{equation*}
	\nabla_{ t}F\left(\frac{t_n}{2},t_n,\frac{t_n}{2}\right)\cdot
	\begin{pmatrix}
		1 \\ 2 \\ 1
	\end{pmatrix}
	= 0.
\end{equation*}	
\end{lemma}

\begin{proof}
We will prove Lemma \ref{lem:magic} by direct calculation. Let us 
begin by estimating $F_\omega\left(\frac{t_n}{2},t_n,\frac{t_n}{2}\right)$. 
Since $F(t,t_n,\frac{t_n}{2})= 0$,  \eqref{eq:fomega}
tells us that
\begin{equation}\label{eq:easy1}
	F_\omega\left(\frac{t_n}{2},t_n,\frac{t_n}{2}\right) 
	= - G^\alpha(\Grad \vvt, \Curl \Lambda^{-\beta}\vvt).
\end{equation}
Similarly, we see that \eqref{eq:fplane} yields
\begin{equation}\label{eq:easy2}
	\begin{split}
		F_\tau = \frac{1}{2}G^\alpha(\Grad \vvt, \Curl \Lambda^{-\beta}\vvt).
	\end{split}
\end{equation}

It remains to estimate $F_t\left(\frac{t_n}{2},t_n,\frac{t_n}{2}\right)$. 
However, as $F(t,t_n,\frac{t_n}{2}) = 0$, 
we must have $F_t\left(\frac{t_n}{2},t_n,\frac{t_n}{2}\right)= 0$. 
This, together with \eqref{eq:easy1}
and \eqref{eq:easy2}, concludes the proof.
\end{proof}

Using the previous lemma, we can now prove that 
the error produced along the diagonal $(t/2,t,t/2)$
is second order in $\Delta t$. 
\begin{lemma}\label{lem:spoton}
	Let $\vvt$ be the Strang splitting solution in 
	the sense of Definition \ref{def:strang} and \eqref{eq:strang-2}.
	Then, for any $k$ such that $\|\vt_0\|_{H^{k+3 \alpha}} \leq C$,
	\begin{equation}
		\left\|F\left(\frac{t}{2},t,\frac{t}{2}\right)\right\|_{H^k} \leq C(\Delta t)^2.
	\end{equation}
\end{lemma}
\begin{proof}
Since $F\left(\frac{t_n}{2},t_n,\frac{t_n}{2}\right) = 0$, a Taylor 
expansion provides the identity
\begin{equation}\label{eq:jeezus}
	\begin{split}
		F\left(\frac{t}{2},t,\frac{t}{2}\right)  &= 
		\nabla_{t}F\left(\frac{t_n}{2},t_n,\frac{t_n}{2}\right)\cdot
		\begin{pmatrix}
			1 \\ 2 \\ 1
		\end{pmatrix}\left(\frac{t}{2}-\frac{t_n}{2}\right) \\
		&\qquad\qquad +\frac{1}{2}\int_{t_n/2}^{t/2}
			\begin{pmatrix}
				1 \\ 2 \\ 1
			\end{pmatrix}^T
			\nabla_{t}^2 F \left(\frac{s}{2},s,\frac{s}{2}\right)
		\begin{pmatrix}
					1 \\ 2 \\ 1
		\end{pmatrix}\left(\frac{s}{2}-\frac{t_n}{2}\right)~ds\\
		&= \frac{1}{2}\int_{t_n/2}^{t/2}
			\begin{pmatrix}
				1 \\ 2 \\ 1
			\end{pmatrix}^T
			\nabla_{t}^2 F \left(\frac{s}{2},s,\frac{s}{2}\right)
		\begin{pmatrix}
					1 \\ 2 \\ 1
		\end{pmatrix}\left(\frac{s}{2}-\frac{t_n}{2}\right)~ds,
	\end{split}
\end{equation}
where the last equality is an application of Lemma \ref{lem:magic}.
By taking the $H^k$ norm on both 
sides of \eqref{eq:jeezus} and applying the 
previous lemma, we gather
\begin{equation*}
	\left\|F\left(\frac{t}{2},t,\frac{t}{2}\right)\right\|_{H^k} \leq C\Delta t^2,
\end{equation*}
which concludes the proof.
\end{proof}

\subsection{Proof of Theorem \ref{thm:rate}}

We have now gathered all the ingredients needed 
to prove the sought after second-order 
error estimate. 

Performing the same calculations as in \eqref{eq:w1}, \eqref{eq:w2}, \eqref{eq:w3},
and then applying Lemma \ref{lem:spoton}, yields
\begin{equation*}
	\frac{1}{2}\partial_t\|e(t)\|_{H^{k-3\alpha}} \leq 
	C(T)\left(\Delta t^2 
	+ \|e\|_{H^{k-3\alpha}}\right), 
	\quad t \in [0,T].
\end{equation*}
Since $e(0)=0$, an application of the Gronwall inequality to the previous inequality gives
\begin{equation*}
	\|e(t)\|_{H^{k-3\alpha}} \leq t\, \Delta t^2C, \quad t \in [0,T],
\end{equation*}
which concludes the proof of Theorem \ref{thm:rate}. \qed

\appendix
\section{Proof of Lemmas \ref{lem:nonlineardiv} and \ref{lem:lineardiv}}
The purpose of this appendix is to provide proofs
of Lemmas \ref{lem:nonlineardiv} and \ref{lem:lineardiv}. 
Both lemmas have been crucial to our convergence analysis. 
In particular, Lemmas \ref{lem:higher}, 
\ref{lem:w}, \ref{lem:F2}, and \ref{lem:tempv}, all rely on 
their validity.
\begin{lemma}\label{lem:nonlineardivA}
Let $k\geq 6$ be an integer. Then
\begin{equation}
	\sum_{s=0}^k\left|\int_{\R^N} \Grad^s(\Grad f\cdot \Curl \Lambda^{-\beta}f):\Grad^s f~dx\right| \leq C\|f\|_{H^{k-2}}\|f\|_{H^k}^2,
\end{equation}	
for all $f \in H^k$.
\end{lemma}

\begin{proof}
Let us confine to the three-dimensional case ($N=3$) 
as the other cases are almost identical. 

By applying the  Leibniz 
rule (with multiindex notation $\alpha = (\alpha_1, \alpha_2, \alpha_3)$), we obtain the following expression 
\begin{align}\label{eq:I1start}
		&\sum_{|\alpha| = s}\int_\Om \Grad^\alpha(\Grad f \cdot \Curl \Lambda^{-\beta}f)\Grad^\alpha f~dx \nonumber\\
		& =  \sum_{|\alpha| = s}\sum_{i_1=0}^{\alpha_1}\sum_{i_2=0}^{\alpha_2}\sum_{i_3=0}^{\alpha_3}
		\left(\alpha_1 \atop i_1\right)\left(\alpha_2 \atop i_2\right)\left(\alpha_3 \atop i_3\right) \\
		& \times\!\! \int_{\R^N} \left(\Grad \frac{\partial^{i_1 + i_2 + i_3}f}{\partial x^{i_1}\partial y^{i_2}\partial z^{i_3}} 
			\cdot \Curl \Lambda^{-\beta}\left(\frac{\partial^{s-i_1 - i_2 - i_3}f}{\partial x^{\alpha_1 - i_1}\partial y^{\alpha_2 - i_2}\partial z^{\alpha_3 - i_3}}\right) \right)
			\frac{\partial^s f}{\partial x^{\alpha_1}\partial y^{\alpha_2}\partial z^{\alpha_3}}~dx. \nonumber
\end{align}
Let us now consider four separate cases of $i_1+i_2+i_3$ in the above quadruple sum.

(i) If $i_1+i_2+ i_3= s$ (i.e $(\alpha_1, \alpha_2, \alpha_3)=(i_1,i_2,i_3)$ ),  the above term can be rewritten as follows
\begin{equation*}
	\begin{split}
		&\int_{\R^N} \left(\Grad \frac{\partial^s f}{\partial x^{\alpha_1}\partial y^{\alpha_2}\partial z^{\alpha_3}} \cdot \Curl \Lambda^{-\beta}f\right)
		\frac{\partial^s f}{\partial x^{\alpha_1}\partial y^{\alpha_2}\partial z^{\alpha_3}}~dx \\
		&\qquad =
		\int_{\R^N} \frac{1}{2}\Grad \left|\frac{\partial^s f}{\partial x^{\alpha_1}\partial y^{\alpha_2}\partial z^{\alpha_3}}\right|^2
		\cdot \Curl \Lambda^{-\beta}f~dx = 0.
	\end{split}
\end{equation*}

(ii) If $2 \leq i_1 + i_2 + i_3  \leq k-3$
\begin{equation*}
	\begin{split}
		&\int_{\R^N} \left(\Grad \frac{\partial^{i_1 + i_2 + i_3}f}{\partial x^{i_1}\partial y^{i_2}\partial z^{i_3}} 
			\cdot \Curl \Lambda^{-\beta}\left(\frac{\partial^{s-i_1 - i_2 - i_3}f}{\partial x^{\alpha_1 - i_1}\partial y^{\alpha_2 - i_2}\partial z^{\alpha_3 - i_3}}\right) \right)
			\frac{\partial^s f}{\partial x^{\alpha_1}\partial y^{\alpha_2}\partial z^{\alpha_3}}~dx \\
		&\leq\! \left\|\Grad \frac{\partial^{i_1 + i_2 + i_3}f}{\partial x^{i_1}\partial y^{i_2}\partial z^{i_3}} \right\|_{L^\infty}
			 \! \left\|\Curl \Lambda^{-\beta}\left(\frac{\partial^{s-i_1 - i_2 - i_3}f}{\partial x^{\alpha_1 - i_1}\partial y^{\alpha_2 - i_2}\partial z^{\alpha_3 - i_3}}\right)\right\|_{L^2}
			\!\left\|\frac{\partial^s f}{\partial x^{\alpha_1}\partial y^{\alpha_2}\partial z^{\alpha_3}}\right\|_{L^2}
		\\ &\leq
			  C\|f\|_{H^{k}}\|f\|_{H^{k-2}}\|f\|_{H^{k}}.
	\end{split}
\end{equation*}

(iii) If $k \leq i_1 + i_2 + i_3 \leq k-1$
\begin{equation*}
	\begin{split}
		&\int_{\R^N} \left(\Grad \frac{\partial^{i_1 + i_2 + i_3}f}{\partial x^{i_1}\partial y^{i_2}\partial z^{i_3}} 
			\cdot \Curl \Lambda^{-\beta}\left(\frac{\partial^{s-i_1 - i_2 - i_3}f}{\partial x^{\alpha_1 - i_1}\partial y^{\alpha_2 - i_2}\partial z^{\alpha_3 - i_3}}\right) \right)
			\frac{\partial^s f}{\partial x^{\alpha_1}\partial y^{\alpha_2}\partial z^{\alpha_3}}~dx \\
		&\leq\! \left\|\Grad \frac{\partial^{i_1 + i_2 + i_3}f}{\partial x^{i_1}\partial y^{i_2}\partial z^{i_3}}\right\|_{L^2}
		\!\left\|\Curl \Lambda^{-\beta}\left(\frac{\partial^{s-i_1 - i_2 - i_3}f}{\partial x^{\alpha_1 - i_1}\partial y^{\alpha_2 - i_2}\partial z^{\alpha_3 - i_3}}\right)\right\|_{L^\infty}
		\!\left\|\frac{\partial^s f}{\partial x^{\alpha_1}\partial y^{\alpha_2}\partial z^{\alpha_3}}\right\|_{L^2} \\
		&\leq C\|f\|_{H^{k}}\|f\|_{H^4}\|f\|_{H^{k}} \leq C\|f\|_{H^{k-2}}\|f\|_{H^{k}}^2.
	\end{split}
\end{equation*}

(iv) If $0 \leq i_1 + i_2 + i_3 \leq 1$
\begin{equation*}
	\begin{split}
		&\int_{\R^N} \left(\Grad \frac{\partial^{i_1 + i_2 + i_3}f}{\partial x^{i_1}\partial y^{i_2}\partial z^{i_3}} 
			\cdot \Curl \Lambda^{-\beta}\left(\frac{\partial^{s-i_1 - i_2 - i_3}f}{\partial x^{\alpha_1 - i_1}\partial y^{\alpha_2 - i_2}\partial z^{\alpha_3 - i_3}}\right) \right)
			\frac{\partial^s f}{\partial x^{\alpha_1}\partial y^{\alpha_2}\partial z^{\alpha_3}}~dx \\
			&\leq\! \left\|\Grad \frac{\partial^{i_1 + i_2 + i_3}f}{\partial x^{i_1}\partial y^{i_2}\partial z^{i_3}}\right\|_{L^\infty}
			\!\left\|\Curl \Lambda^{-\beta}\left(\frac{\partial^{s-i_1 - i_2 - i_3}f}{\partial x^{\alpha_1 - i_1}\partial y^{\alpha_2 - i_2}\partial z^{\alpha_3 - i_3}}\right)\right\|_{L^2}
			\!\left\|\frac{\partial^s f}{\partial x^{\alpha_1}\partial y^{\alpha_2}\partial z^{\alpha_3}}\right\|_{L^2} \\
			&\leq C\|f\|_{H^4}\|f\|_{H^{k}}^2 \leq C\|f\|_{H^{k-2}}\|f\|_{H^{k}}^2.
	\end{split}
\end{equation*}

Hence, by applying (i)-(iv) in \eqref{eq:I1start}, we see that
\begin{equation}\label{eq:I1bound}
	\sum_{|\alpha| = s}\int_\Om \Grad^\alpha(\Grad f \cdot \Curl \Lambda^{-\beta}f)\Grad^\alpha f~dx  \leq C\|f\|_{H^{k-2}}\|f\|_{H^{k}}^2,
\end{equation}
which concludes the proof.
\end{proof}

\begin{lemma}\label{lem:lineardivA}
Let $k\geq 4$ be an integer. The following estimates hold 
	\begin{align}\label{eq:ldiv1A}
		\sum_{s=0}^k\left|\int_{\R^N} \Grad^s \left(\Grad f \cdot \Curl \Delta^{-\beta} g  \right):\Grad^s f~dx\right| &\leq C\|g\|_{H^k}\|f\|_{H^k}^2, \quad f,g\in H^k, \\
\label{eq:ldiv2A}
			\sum_{s=0}^k\left|\int_{\R^N}\Grad^s\left(\Grad g \cdot \Curl \Lambda^{-\beta}f\right):\Grad^s f~dx\right| &\leq C\|g\|_{H^{k+1}}\|f\|_{H^k}^2,\quad f\in H^k, \, g\in H^{k+1}.
	\end{align}
\end{lemma}

\begin{proof}
	The proof of \eqref{eq:ldiv1} is easily obtained by the calculations
	of the previous proof. To prove \eqref{eq:ldiv2}, it is only step (i)
 	of the previous proof which is no longer true. However, 
	this is also the reason for the $k+1$ on $g$. That is, step (i)
	is replaced by a simpler H\"older inequality.
	
\end{proof}

\end{document}